\documentclass[a4paper,12pt,reqno]{amsart}
\usepackage{amsmath}
\usepackage{amssymb}
\usepackage{amsthm}

\usepackage{amscd}
\usepackage{amsfonts}
\usepackage{setspace}
\usepackage{version}

\title[High moments of random multiplicative functions]{Moments of random multiplicative functions, II: High moments}
\author{Adam J Harper}
\address{Mathematics Institute, Zeeman Building, University of Warwick, Coventry CV4 7AL, England}
\email{A.Harper@warwick.ac.uk}
\date{11th April 2018}
\thanks{When this work was started, the author was supported by a research fellowship at Jesus College, Cambridge. The work continued while the author was in residence at the Mathematical Sciences Research Institute in Berkeley, California (supported by the National Science Foundation under Grant No. DMS-1440140), during the Spring 2017 semester.}

\numberwithin{equation}{section}
\theoremstyle{plain}

\onehalfspace
\newcommand{\N}{\mathbb{N}}
\newcommand{\R}{\mathbb{R}}
\newcommand{\E}{\mathbb{E}}
\newcommand{\p}{\mathbb{P}}
\newcommand{\Z}{\mathbb{Z}}

\newtheorem{thmhigh1}{Theorem}
\newtheorem{thmhigh2}[thmhigh1]{Theorem}
\newtheorem{cor1}{Corollary}

\newtheorem{epres1}{Euler Product Result}
\newtheorem{epres2}[epres1]{Euler Product Result}

\newtheorem{probres1}{Probability Result}
\newtheorem{probres2}[probres1]{Probability Result}
\newtheorem{probres3}[probres1]{Probability Result}

\newtheorem{harman1}{Harmonic Analysis Result}
\newtheorem{numth1}{Number Theory Result}

\newtheorem{prop1}{Proposition}
\newtheorem{prop2}[prop1]{Proposition}
\newtheorem{prop3}[prop1]{Proposition}
\newtheorem{prop4}[prop1]{Proposition}

\newtheorem{keyprop1}{Key Proposition}
\newtheorem{keyprop2}[keyprop1]{Key Proposition}

\begin{document}

\maketitle

\begin{abstract}
We determine the order of magnitude of $\E|\sum_{n \leq x} f(n)|^{2q}$ up to factors of size $e^{O(q^2)}$, where $f(n)$ is a Steinhaus or Rademacher random multiplicative function, for all real $1 \leq q \leq \frac{c\log x}{\log\log x}$.

In the Steinhaus case, we show that $\E|\sum_{n \leq x} f(n)|^{2q} = e^{O(q^2)} x^q (\frac{\log x}{q\log(2q)})^{(q-1)^2}$ on this whole range. In the Rademacher case, we find a transition in the behaviour of the moments when $q \approx (1+\sqrt{5})/2$, where the size starts to be dominated by ``orthogonal'' rather than ``unitary'' behaviour. We also deduce some consequences for the large deviations of $\sum_{n \leq x} f(n)$.

The proofs use various tools, including hypercontractive inequalities, to connect $\E|\sum_{n \leq x} f(n)|^{2q}$ with the $q$-th moment of an Euler product integral. When $q$ is large, it is then fairly easy to analyse this integral. When $q$ is close to 1 the analysis seems to require subtler arguments, including Doob's $L^p$ maximal inequality for martingales.
\end{abstract}

\section{Introduction}
In this sequence of papers, we are interested in the moments $\E|\sum_{n \leq x} f(n)|^{2q}$ of random multiplicative functions $f(n)$.

We consider two different models for $f(n)$, a {\em Steinhaus random multiplicative function} and a {\em Rademacher random multiplicative function}. We obtain a Steinhaus random multiplicative function by letting $(f(p))_{p \; \text{prime}}$ be a sequence of independent Steinhaus random variables (i.e. distributed uniformly on the unit circle $\{|z|=1\}$), and then setting $f(n) := \prod_{p^{a} || n} f(p)^{a}$ for all natural numbers $n$, where $p^a || n$ means that $p^a$ is the highest power of the prime $p$ that divides $n$. We obtain a Rademacher random multiplicative function by letting $(f(p))_{p \; \text{prime}}$ be independent Rademacher random variables (i.e. taking values $\pm 1$ with probability $1/2$ each), and then setting $f(n) := \prod_{p |n} f(p)$ for all squarefree $n$, and $f(n) = 0$ when $n$ is not squarefree.

Random multiplicative functions have attracted quite a lot of attention as models for functions of number theoretic interest: for example, Rademacher random multiplicative functions were introduced by Wintner~\cite{wintner} as a model for the M\"{o}bius function $\mu(n)$. There are also probabilistic and analytic motivations for studying them, see Saksman and Seip's open problems paper~\cite{saksmanseip}, for example. The introduction to the previous paper~\cite{harperrmflow} in this sequence contains a more extensive discussion of some of these connections.

\vspace{12pt}
Harper~\cite{harperrmflow} showed that for Steinhaus or Rademacher random multiplicative $f(n)$, for all large $x$ we have
$$ \E|\sum_{n \leq x} f(n)|^{2q} \asymp \left( \frac{x}{1 + (1-q)\sqrt{\log\log x}} \right)^q \;\;\;\;\; \forall \; 0 \leq q \leq 1 . $$
In particular, taking $q = 1/2$ this implies that $\E|\sum_{n \leq x} f(n)| \asymp \frac{\sqrt{x}}{(\log\log x)^{1/4}}$, which proved a conjecture of Helson~\cite{helson} that the first absolute moment should be $o(\sqrt{x})$.

Our goal here is to investigate the case where $q \geq 1$. When $q \in \N$ is fixed, one can expand the $2q$-th power and reduce the calculation of $\E|\sum_{n \leq x} f(n)|^{2q}$ to a number theoretic counting problem. For example, in the Steinhaus case one has
$$ \E|\sum_{n \leq x} f(n)|^{2q} = \#\{n_1, ..., n_{2q} \leq x : \prod_{i=1}^{q} n_i = \prod_{i=q+1}^{2q} n_i\} . $$
Starting from this, one can obtain an asymptotic for the moment as $x \rightarrow \infty$, which was carried out by Harper, Nikeghbali and Radziwi{\l}{\l}~\cite{hnr}, and also independently by Heap and Lindqvist~\cite{heaplindqvist}, and (in the Steinhaus case) in unpublished work of Granville and Soundararajan. The result is that, for fixed $q \in \N$ and Steinhaus random multiplicative $f(n)$, one has
\begin{equation}\label{steinhausasymp}
\E|\sum_{n \leq x} f(n)|^{2q} \sim C_{\text{St}}(q) x^{q} \log^{(q-1)^2}x \;\;\; \text{as} \; x \rightarrow \infty ,
\end{equation}
where the constant $C_{\text{St}}(q)$ satisfies $C_{\text{St}}(q) = e^{-q^2 \log q - q^{2}\log\log q + O(q^2)}$ for large $q$. For Rademacher random multiplicative $f(n)$, when $q=1$ we have that $\E|\sum_{n \leq x} f(n)|^{2} = \sum_{n \leq x, \; n \; \text{squarefree}} 1 \sim (6/\pi^2)x$, and for fixed integer $q \geq 2$ we have
$$ \E|\sum_{n \leq x} f(n)|^{2q} \sim C_{\text{Rad}}(q) x^{q} \log^{q(2q-3)}x \;\;\; \text{as} \; x \rightarrow \infty , $$
where the constant $C_{\text{Rad}}(q)$ satisfies $C_{\text{Rad}}(q) = e^{-2q^2 \log q - 2q^{2}\log\log q + O(q^2)}$ for large $q$. As described in \cite{hnr,heaplindqvist}, we actually have much more precise information about the constants $C_{\text{St}}(q), C_{\text{Rad}}(q)$ (for example they factor into explicit ``arithmetic'' and ``geometric'' parts), but this will not be important for our purposes here.

We would like to have information about $\E|\sum_{n \leq x} f(n)|^{2q}$ when $q \geq 1$ is not necessarily integral, and that allows $q$ to vary as a function of $x$ rather than being fixed.

Regarding uniformity in $q$, Theorem 4.1 of Granville and Soundararajan~\cite{gransoundlcs} implies that for Steinhaus random multiplicative $f(n)$, and uniformly for all large $x$ and integers $q \geq 1$ such that $q^{eq} \leq x$, we have
$$ e^{-q^{2}\log q - q^{2}\log\log(2q)} (\log(\frac{\log x}{q\log 2q}))^{-O(q^2)} \leq \frac{\E|\sum_{n \leq x} f(n)|^{2q}}{x^{q} \log^{(q-1)^2}x} \leq e^{-q^{2}\log q + O(q^2)} . $$
This range of $q$ is essentially the largest on which one could expect a result of a similar shape to \eqref{steinhausasymp}. Indeed, if $q \geq A \frac{\log x}{\log\log x}$ for some $A \geq 1$ (say) then we have $e^{-q^2 \log q - q^{2}\log\log(2q)} x^q \log^{(q-1)^2}x \leq ((1+o(1))A)^{-q^2} x^q$, which becomes incompatible with the lower bound $\E |\sum_{n \leq x} f(n)|^{2q} \geq (\E |\sum_{n \leq x} f(n)|^{2})^q = \lfloor x \rfloor^q$ coming from H\"{o}lder's inequality\footnote{In this paper we are not particularly concerned with the case where $q \geq \frac{\log x}{\log\log x}$, but for completeness we make a few indicative remarks. Section 6 of Granville and Soundararajan~\cite{gransoundlcs} contains various results on this range of $q$. Setting $v = \log(2q(\log q)/\log x) \gg 1$, and redoing the calculations in section \ref{secsteineasy} with the Rankin shift $1 + \frac{q}{\log x}$ replaced by $1 + \frac{v}{\log(q\log x)}$ and with $q^2$-smooth numbers replaced by $q\log x$-smooth numbers, one can show that $\E |\sum_{n \leq x} f(n)|^{2q} \leq x^{q(1 + \frac{v}{\log(q\log x)} + o(1))}$ uniformly for $q \geq \frac{\log x}{\log\log x}$.  In particular, if $q = \log^{1+a}x$ for any fixed $a \geq 0$ then we have $\E |\sum_{n \leq x} f(n)|^{2q} \leq x^{q(1 + \frac{a}{a+2} + o(1))}$. By only considering the contribution to the expectation from the event that $f(p)$ is very close to 1 for all primes $p \leq \frac{q\log x}{\log\log x} = \frac{\log^{2+a}x}{\log\log x}$, one can obtain a comparable lower bound for $\E |\sum_{n \leq x} f(n)|^{2q}$ (as in Corollary 6.3 of Granville and Soundararajan~\cite{gransoundlcs}).}. But the bounds are imperfect, as the upper bound doesn't include the factor $e^{-q^{2}\log\log(2q)}$ that we expect to appear, and the lower bound features the extraneous factor $(\log(\frac{\log x}{q\log 2q}))^{-O(q^2)}$. They also remain restricted to {\em integer} $q$. There are various other results in the literature that study the Steinhaus moments $\E|\sum_{n \leq x} f(n)|^{2q}$, and variants of them, for integer $q$, especially for small integers where one can try to obtain lower order terms in the known asymptotics. See e.g. the preprint of Shi and Weber~\cite{shiweber}, and the references cited there. However, the author is not aware of any work giving sharp moment bounds for non-integer $q$, nor improving the dependence on $q$ in Granville and Soundararajan's~\cite{gransoundlcs} bounds for the large integer case.

We shall prove the following uniform estimate for all {\em real} $q$.
\begin{thmhigh1}
There exists a small absolute constant $c > 0$ such that the following is true. If $f(n)$ is a Steinhaus random multiplicative function, then uniformly for all large $x$ and real $1 \leq q \leq \frac{c\log x}{\log\log x}$ we have
$$ \E |\sum_{n \leq x} f(n)|^{2q} = e^{-q^{2}\log q - q^{2}\log\log(2q) + O(q^2)} x^{q} \log^{(q-1)^2}x . $$
\end{thmhigh1}

To avoid any confusion, we restate this first result more explicitly: on the stated range of $q$ and $x$, we always have
$$ e^{-q^{2}\log q - q^{2}\log\log(2q) - Cq^2} \leq \frac{\E |\sum_{n \leq x} f(n)|^{2q}}{x^{q} \log^{(q-1)^2}x} \leq e^{-q^{2}\log q - q^{2}\log\log(2q) + Cq^2} , $$
for a certain absolute constant $C$. We do {\em not} know how to prove an asymptotic like \eqref{steinhausasymp} when $q$ is not a fixed natural number.

In the Rademacher case, even conjecturally the behaviour of $\E|\sum_{n \leq x} f(n)|^{2q}$ is perhaps not obvious. On a wide range of real $q \geq 2$, we might expect that $\E |\sum_{n \leq x} f(n)|^{2q} = e^{-2q^{2}\log q - 2q^{2}\log\log(2q) + O(q^2)} x^{q} \log^{q(2q-3)}x$ as in the known asymptotics. But this certainly cannot be the answer for all $1 \leq q \leq 2$, since on some of that range the exponent $q(2q-3)$ of the logarithm would be negative. (And, by H\"{o}lder's inequality, we must at least have $\E |\sum_{n \leq x} f(n)|^{2q} \geq (\E |\sum_{n \leq x} f(n)|^{2})^q \gg x^q$.)
\begin{thmhigh2}
Let $q_0 = (1 + \sqrt{5})/2 \approx 1.618$. There exists a small absolute constant $c > 0$ such that the following is true. If $f(n)$ is a Rademacher random multiplicative function, then uniformly for all large $x$ and real $1 \leq q \leq \frac{c\log x}{\log\log x}$ we have
$$ \E |\sum_{n \leq x} f(n)|^{2q} = e^{-2q^{2}\log q - 2q^{2}\log\log(2q) + O(q^2)} (1 + \min\{\log\log x , \frac{1}{|q - q_0 |}\}) x^{q} \log^{\max\{(q-1)^2, q(2q-3)\}}x . $$
\end{thmhigh2}

With hindsight, the exponent of $\log x$ we obtain in Theorem 2 is perhaps quite natural, since one doesn't expect slower growth in the Rademacher than the Steinhaus case (where there is ``more room'' for the complex valued random variables to cancel), and we expect $q(2q-3)$ to be the correct exponent eventually. Notice that the golden ratio $q_0$ is the value at which $q(2q-3)$ becomes larger than $(q-1)^2$. But the additional factor $\min\{\log\log x , \frac{1}{|q - q_0 |}\}$ that appears for $q$ close to $q_0$ seems genuinely unexpected, and hard to understand except through an inspection of the proof of the theorem.

\vspace{12pt}
Next we shall discuss the proofs. Once $q$ is moderately large, namely when $q \geq \log\log x$, we can prove the upper bounds in Theorems 1 and 2 by fairly simple arguments. See section \ref{seceasycases}. This is because, for such $q$, terms like $\log^{O(q)}x$ can be absorbed into the factor $e^{O(q^2)}$ in our theorems, so we can afford to use simple techniques that are a bit wasteful (e.g. involving H\"{o}lder's inequality to reduce to the case of integer $q$) to reduce matters to a counting problem. Then Rankin's trick is almost sufficient to perform the relevant counts. To obtain the terms $e^{- q^{2}\log\log(2q)}$ and $e^{- 2q^{2}\log\log(2q)}$ in the theorems, we use Rankin's trick along with a slightly more careful treatment of small prime factors.

Our main work is to prove Theorems 1 and 2 for $1 \leq q \leq \log\log x$, and also the lower bounds for larger $q$. Let $F(s) = \sum_{\substack{n=1, \\ p|n \Rightarrow p \leq x}}^{\infty} \frac{f(n)}{n^s}$ denote the Dirichlet series corresponding to $f(n)$, on $x$-smooth numbers (i.e. numbers with all their prime factors $\leq x$). We can also write $F(s)$ as an Euler product, namely $F(s) = \prod_{p \leq x} (1 - \frac{f(p)}{p^{s}})^{-1}$ in the Steinhaus case and $F(s) = \prod_{p \leq x} (1 + \frac{f(p)}{p^{s}})$ in the Rademacher case. In the author's treatment~\cite{harperrmflow} of low moments, the first step was to show (roughly) that $\E|\sum_{n \leq x} f(n)|^{2q} \approx x^{q} \E\left(\frac{1}{\log x} \int_{-1/2}^{1/2} |F(1/2+it)|^2 dt \right)^{q}$ when $2/3 \leq q \leq 1$. Similarly, our first step here is to show that
\begin{equation}\label{introprod}
\E|\sum_{n \leq x} f(n)|^{2q} \approx e^{O(q^2)} x^{q} \E\Biggl(\frac{1}{\log x} \int_{-1/2}^{1/2} |F(1/2+ \frac{q}{\log x} + it)|^2 dt \Biggr)^{q} .
\end{equation}
Note the shift by $q/\log x$ in the integral, which is analogous to the use of Rankin's trick in our elementary upper bound argument for $q \geq \log\log x$. The basic strategy for proving something like \eqref{introprod} is the same as in \cite{harperrmflow}, namely conditioning on the behaviour of $f(n)$ on smaller primes; using fairly standard moment inequalities, like Khintchine's inequality, to show the conditional expectation behaves like a power of a mean square average; and using Parseval's identity to relate the mean square to an integral average of the Euler product. In \cite{harperrmflow} one could bound terms by using H\"{o}lder's inequality to pass to the second moment, whereas here we need suitable rough bounds for high moments. These are supplied by a pair of hypercontractive inequalities, see Probability Result 1 in section \ref{secprelimresults}. Applying the hypercontractive inequalities introduces various divisor functions $d_{\lceil q \rceil}(n), d_{2\lceil q \rceil - 1}(n)$ into our calculations, requiring a bit more number theoretic work as compared with the low moments argument of \cite{harperrmflow}. We refer to the beginning of section \ref{seceulerprod} for a rigorous formulation of \eqref{introprod}, and more technical comparison of this part of the argument with the low moments case~\cite{harperrmflow}.

Next, we observe that the right hand side of \eqref{introprod} is
\begin{equation}\label{introproddisc}
\approx e^{O(q^2)} x^{q} \E\Biggl(\frac{1}{\log^{2}x} \sum_{|n| \leq (\log x)/2} |F(1/2+ \frac{q}{\log x} + i\frac{n}{\log x})|^2 \Biggr)^{q} ,
\end{equation}
since heuristically the value of $|F(1/2+ \frac{q}{\log x} + it)|$ doesn't change much on $t$ intervals of length $1/\log x$. One can obtain rigorous statements of this kind using H\"{o}lder's inequality in the upper bound arguments, and Jensen's inequality in the lower bound arguments, see sections \ref{secmainupper} and \ref{secmainlower}. Now we can see heuristically why Theorems 1 and 2 might hold. In the Steinhaus case, the Euler product $F(s)$ behaves on average like an $L$-function from a {\em unitary} family, and then since we have $q \geq 1$ (and very differently than in the low moments case~\cite{harperrmflow}) the sum over $n$ essentially gives us $\log x$ independent tries at obtaining a large value of $F(s)$. So the right hand side of \eqref{introproddisc} is $\approx e^{O(q^2)} x^{q} \frac{1}{\log^{2q}x} \log x \E |F(1/2+ \frac{q}{\log x})|^{2q} \approx e^{O(q^2)} x^{q} \frac{1}{\log^{2q-1}x} (\frac{\log x}{q\log(2q)})^{q^2}$, as in Theorem 1. In the Rademacher case, $F(1/2+ \frac{q}{\log x} + it)$ behaves like an $L$-function from an {\em orthogonal} family when $t \approx 0$, and like an $L$-function from a unitary family\footnote{At first glance, one might expect $F(1/2+ \frac{q}{\log x} + it)$ to behave like a {\em symplectic} $L$-function when $t \approx 0$, because averaging over Rademacher $f(n)$ models averaging over quadratic Dirichlet characters. The reason we actually have orthogonal behaviour is because we restrict our sums $\sum_{n \leq x} f(n)$ to squarefree terms. For some other contexts where a transition from orthogonal/symplectic to unitary behaviour arises, as for large $t$ here, see the papers of Florea~\cite{florea}, Keating and Odgers~\cite{keatodg}, and Soundararajan and Young~\cite{soundyoung}, for example.} when $t \approx 1$. Thus, thanks to those $(\log x)/4 \leq |n| \leq (\log x)/2$ (say) we get a contribution $e^{O(q^2)} x^{q} \frac{1}{\log^{2q-1}x} (\frac{\log x}{q\log(2q)})^{q^2}$ to the right hand side of \eqref{introproddisc}, and thanks to the $n=0$ term we get a contribution $\approx e^{O(q^2)} x^{q} \frac{1}{\log^{2q}x} \E |F(1/2+ \frac{q}{\log x})|^{2q} \approx e^{O(q^2)} x^{q} \frac{1}{\log^{2q}x} (\frac{\log x}{q\log(2q)})^{2q^2 - q}$. The factor $(1 + \min\{\log\log x , \frac{1}{|q - q_0 |}\})$ in Theorem 2 arises because of the contribution from intermediate values of $n$.

To prove the lower bounds in Theorems 1 and 2 rigorously, as we do in section \ref{secmainlower}, roughly speaking it suffices to note that \eqref{introproddisc} is $\geq e^{O(q^2)} x^{q} \frac{1}{\log^{2q}x} \E \sum_{|n| \leq (\log x)/2} |F(1/2+ \frac{q}{\log x} + i\frac{n}{\log x})|^{2q}$, and then compute $\E |F(1/2+ \frac{q}{\log x} + i\frac{n}{\log x})|^{2q}$. In practice the details are slightly more complicated because the precise version of \eqref{introproddisc} involves some other terms, including subtracted error terms that must be upper bounded. However, we can obtain suitable upper bounds from our main section \ref{secmainupper} argument for proving the upper bounds in Theorems 1 and 2.

To prove those upper bounds rigorously, we need to capture the fact that typically there will only be a few large terms in the sum over $n$ in \eqref{introproddisc}. When $q \geq 2$, a careful application of H\"{o}lder's inequality lets us bound \eqref{introproddisc} by estimating terms of the form $\E |F(1/2+ \frac{q}{\log x} + i\frac{n}{\log x})|^2  |F(1/2+ \frac{q}{\log x} + i\frac{m}{\log x})|^{2(q-1)}$. These decrease in size quite rapidly as $|m-n|$ becomes large (and, in the Rademacher case, as $|m|, |n|$ become large), because the parts of the two Euler products over primes $> x^{1/|m-n|}$ become decorrelated rather than reinforcing one another. This indeed says that one doesn't expect large contributions from many different $m,n$. When $1 < q < 2$, such a direct argument doesn't seem to succeed, so we need a more subtle approach. The rough idea is to treat parts of the Euler products over ``small'' and ``large'' primes differently, so after a (different) careful application of H\"{o}lder's inequality, one is led to expectations where different parts of the Euler product appear to different exponents, to maximise the decorrelation we capture. The most difficult situation is where $q$ is very close to 1 (i.e. $q=1 + o(1)$ as $x \rightarrow \infty$). To handle this without picking up any terms that blow up as $q$ approaches 1, we use a martingale maximal inequality (see Probability Result 3 in section \ref{secprelimresults}) that essentially lets us maximise over several different splittings of the Euler product simultaneously.

\vspace{12pt}
As just described, we go to quite a lot of trouble to prove Theorems 1 and 2 when $q$ is just a little larger than 1. It is satisfying to have a uniform result (and a method capable of proving one), but in addition this range of $q$ turns out to be relevant for deducing the following corollary.

\begin{cor1}
Let $x$ be large, and let $f(n)$ be a Steinhaus or Rademacher random multiplicative function. For all $2 \leq \lambda \leq \sqrt{\log x}$, say, we have
$$ \p(|\sum_{n \leq x} f(n)| \geq \lambda \sqrt{x}) \ll \frac{1}{\lambda^2} e^{-(\log^{2}\lambda)/\log\log x} . $$
\end{cor1}

\begin{proof}[Proof of Corollary 1]
For any $1 \leq q \leq 3/2$, say, Theorems 1 and 2 imply that
$$ \p(|\sum_{n \leq x} f(n)| \geq \lambda \sqrt{x}) \leq \frac{\E|\sum_{n \leq x} f(n)|^{2q}}{(\lambda \sqrt{x})^{2q}} \ll \frac{\log^{(q-1)^2}x}{\lambda^{2q}} = \frac{1}{\lambda^{2}} e^{(q-1)^{2}\log\log x - 2(q-1)\log\lambda} . $$
Calculus implies that the right hand side is minimised if we choose $q-1 = \frac{\log\lambda}{\log\log x}$, and inserting this choice proves Corollary 1. 
\end{proof}

In the paper~\cite{harperrmflow} on low moments, by considering $\E|\sum_{n \leq x} f(n)|^{2q}$ with $q$ a little smaller than 1 the author showed that $\p(|\sum_{n \leq x} f(n)| \geq z\frac{\sqrt{x}}{(\log\log x)^{1/4}}) \ll \frac{\min\{\log z, \sqrt{\log\log x}\}}{z^2}$ for all $z \geq 2$. Corollary 1 is weaker than this when $\lambda \leq e^{\sqrt{\log\log x}}$, but stronger for larger $\lambda$. In \cite{harperrmflow} the author also showed (see Corollary 2 there, and the subsequent discussion) that $\p(|\sum_{n \leq x} f(n)| \geq z\frac{\sqrt{x}}{(\log\log x)^{1/4}}) \gg \frac{e^{-(\log^{2}z)/\log\log x}}{z^{2} (\log\log x)^{O(1)}}$ on a wide range of $z$. Together all these results give a fairly complete description of the tail behaviour of $\sum_{n \leq x} f(n)$, up to factors $(\log\log x)^{O(1)}$.

\vspace{12pt}
We end this introduction with a few remarks on other possible approaches to Theorems 1 and 2, and connections with the wider literature.

The quantity $\frac{1}{\log x} \int_{-1/2}^{1/2} |F(1/2+ \frac{q}{\log x} + it)|^2 dt$ in \eqref{introprod} is closely related to (the total mass of a truncation of) a probabilistic object called {\em critical multiplicative chaos}. This connection is discussed extensively in the introduction to the low moments paper~\cite{harperrmflow}, since in that case the techniques for analysing $\E\left(\frac{1}{\log x} \int_{-1/2}^{1/2} |F(1/2+ \frac{q}{\log x} + it)|^2 dt \right)^{q}$ are heavily motivated by ideas from the multiplicative chaos literature. When $q > 1$ the analogous problem does not seem to have been investigated for critical multiplicative chaos, since the $q$-th moment of the integral will diverge as $x \rightarrow \infty$ and this seems to be all the information that was wanted in that case (where the usual interest is in letting $x \rightarrow \infty$ and obtaining a limiting measure whose properties can be investigated). Theorems 1 and 2 show very different behaviour in the Steinhaus and Rademacher cases when $q$ is large, whereas in the usual problems of multiplicative chaos one finds rather universal behaviour (and indeed the Steinhaus and Rademacher moments are of the same order when $q \leq 1$).

Assuming the Generalised Riemann Hypothesis for Dirichlet $L$-functions,  Munsch~\cite{munsch} proved almost sharp upper bounds for the $2k$-th moment of theta functions $\theta(1,\chi)$ as the character $\chi$ varies over non-principal Dirichlet characters mod $q$, for each fixed $k \in \N$. He did this by writing $\theta(1,\chi)$ as a Perron integral involving the $L$-function $L(s,\chi)$, and then expanding the $2k$-th power and bounding the averages of products $\prod_{j=1}^{2k} |L(1/2+it_j, \chi)|$ that emerge. This is interesting here because for even characters $\chi$, $\theta(1,\chi)$ behaves roughly like $\sum_{n \leq \sqrt{q}} \chi(n)$, which is modelled by the sum $\sum_{n \leq \sqrt{q}} f(n)$ of a Steinhaus random multiplicative function. In our case, using Perron's formula we have
\begin{eqnarray}
\E|\sum_{n \leq x} f(n)|^{2q} & \approx & \E\Biggl| \frac{1}{2\pi} \int_{-\sqrt{x}}^{\sqrt{x}} F(1/2 + \frac{q}{\log x} + it) \frac{x^{1/2+q/\log x + it}}{1/2 + q/\log x + it} dt \Biggr|^{2q} \nonumber \\
& \leq & x^{q} e^{O(q^2)} \E\Biggl( \int_{-\sqrt{x}}^{\sqrt{x}} |F(1/2 + \frac{q}{\log x} + it)| \frac{dt}{|1/2 + q/\log x + it|} \Biggr)^{2q} , \nonumber
\end{eqnarray}
say. We already have asymptotics for $\E|\sum_{n \leq x} f(n)|^{2q}$ for fixed $q \in \N$, but we might hope to get an alternative proof of sharp upper bounds for $q \notin \N$ by using H\"{o}lder's inequality in some way. A direct application, producing a term $|F(1/2 + \frac{q}{\log x} + it)|^{2q}$, cannot give sharp bounds because it doesn't recognise that the size of the expectation will be dominated by the integral of $F(1/2 + \frac{q}{\log x} + it)$ over a very short (random) $t$ interval. To detect this, one could pull out a few (say $d$) copies of the bracket before applying H\"{o}lder's inequality to the remaining ones. This would produce a multiple integral of terms of the form $\E (\prod_{j=1}^{d} |F(1/2 + \frac{q}{\log x} + it_j)|) |F(1/2 + \frac{q}{\log x} + iu)|^{2q-d}$, and the biggest contribution comes when all of the $t_j$ are approximately equal to $u$, so indeed we would capture the localisation of the largest contributions. Based on a few rough calculations, it appears this alternative method can prove sharp upper bounds if we take $d=3$ (we need to pull out enough terms to adequately detect the localisation), and if $q \geq 5$, say. But for smaller $q$ this kind of argument doesn't seem operable to prove sharp bounds, indeed one has already lost too much information in applying the triangle inequality to the Perron integral. Nevertheless, it might permit a relatively straightforward extension of Munsch's~\cite{munsch} results to non-integer $k \geq 5$.

A standard strategy for proving lower bounds is to calculate $\E (\sum_{n \leq x} f(n)) R_{x,q}(f)$ and $\E |R_{x,q}(f)|^{2q/(2q-1)}$, where $R_{x,q}(f)$ is some function that is chosen as a proxy for $(\sum_{n \leq x} f(n))^{2q-1}$ that is easier to understand. Then H\"{o}lder's inequality gives
$$ \E|\sum_{n \leq x} f(n)|^{2q} \geq \frac{|\E (\sum_{n \leq x} f(n)) R_{x,q}(f)|^{2q}}{(\E |R_{x,q}(f)|^{2q/(2q-1)})^{2q-1}} . $$
If we can estimate the expectations in the numerator and denominator, and $R_{x,q}(f)$ is well chosen so that both of them do behave like $\E|\sum_{n \leq x} f(n)|^{2q}$ (up to scaling factors that would cancel out), then one obtains a sharp lower bound for the $2q$-th moment. Munsch and Shparlinski~\cite{munshpar} proved sharp lower bounds for the $2k$-th moments of theta functions $\theta(1,\chi)$, for fixed $k \in \N$, by implementing this strategy with a power of a short character sum chosen as the ``proxy'' object. Our analysis shows that for Rademacher random multiplicative functions, we can imagine heuristically that $|\sum_{n \leq x} f(n)| \approx \frac{\sqrt{x}}{\log x} |F(1/2 + \frac{q}{\log x})|$ (when studying $2q$-th moments with $q > q_0$). Motivated by this, we could try taking $R_{x,q}(f) = |F(1/2+\frac{q}{\log x})|^{2q-1}$, or perhaps a small variant of this where primes smaller than $q^{O(1)}$ are excluded from the Euler product. In the Rademacher case, rough calculations suggest this will indeed yield the lower bound $\E|\sum_{n \leq x} f(n)|^{2q} \geq e^{O(q^2)} x^{q} (\frac{\log x}{q\log(2q)})^{q(2q-3)}$, which is sharp when $q > q_0$. For smaller $q$, and in the Steinhaus case, our analysis suggests taking $R_{x,q}(f) = \sum_{|m| \leq \log x} |F(1/2 + \frac{q}{\log x} + i\frac{m}{\log x})|^{2q-1}$. This choice actually won't quite work, but rough calculations suggest that comparing $\E |\sum_{n \leq x} f(n)|^2 \sum_{|m| \leq \log x} |F(1/2 + \frac{q}{\log x} + i\frac{m}{\log x})|^{2(q-1)}$ and $\E (\sum_{|m| \leq \log x} |F(1/2 + \frac{q}{\log x} + i\frac{m}{\log x})|^{2(q-1)})^{q/(q-1)}$ will yield sharp lower bounds for $\E|\sum_{n \leq x} f(n)|^{2q}$. This does not seem simpler than our original proofs of the lower bounds in Theorems 1 and 2, however.

\subsection{Notation and references}
We will say a number $n$ is $y$-smooth if all prime factors of $n$ are $\leq y$. We will generally use $p$ to denote primes. Unless mentioned otherwise, the letters $c, C$ will denote positive constants, $c$ usually being a small constant and $C$ a large one. We write $f(x) = O(g(x))$ and $f(x) \ll g(x)$, both of which mean that there exists $C$ such that $|f(x)| \leq Cg(x)$, for all $x$. Sometimes this notation will be adorned with a subscript parameter (e.g. $O_{\epsilon}(\cdot)$ and $\ll_{\delta}$), meaning that the implied constant $C$ is allowed to depend on that parameter. We write $f(x) \asymp g(x)$ to mean that $g(x) \ll f(x) \ll g(x)$, in other words that $c g(x) \leq |f(x)| \leq C g(x)$ for some $c,C$, for all $x$.

The books of Gut~\cite{gut} and of Montgomery and Vaughan~\cite{mv} may be consulted as excellent general references for probabilistic and number theoretic background for this paper.

\section{Preliminary results}\label{secprelimresults}

\subsection{Random Euler products}
We begin with some ``two point'' estimates for the expectation of the $2\alpha$-th power of a random Euler product, multiplied by the $2\beta$-th power of an imaginary shift of that product. These estimates, and small variants of them, will be basic tools throughout our work. The calculations are closely related to computations of shifted moments of $L$-functions, as in the papers of Chandee~\cite{chandee} and of Soundararajan and Young~\cite{soundyoung}, for example.

\begin{epres1}
If $f$ is a Steinhaus random multiplicative function, then for any real $\alpha, \beta \geq 0$, any real $100(1 + \max\{\alpha^2 , \beta^2\}) \leq x \leq y$, and any real $\sigma \geq - 1/\log y$ and $t$, we have
\begin{eqnarray}
&& \E \prod_{x < p \leq y} \left|1 - \frac{f(p)}{p^{1/2+\sigma}}\right|^{-2\alpha} \left|1 - \frac{f(p)}{p^{1/2+\sigma + it}}\right|^{-2\beta} \nonumber \\
& = & \exp\{\sum_{x < p \leq y} \frac{\alpha^2 + \beta^2 + 2\alpha\beta \cos(t\log p)}{p^{1 + 2\sigma}} + O(\frac{\max\{\alpha, \beta, \alpha^3 , \beta^3\}}{\sqrt{x} \log x}) \} . \nonumber
\end{eqnarray}

If we also have $\sigma \leq 1/\log y$, then the above is
$$ = e^{O(\max\{\alpha, \beta, \alpha^2 , \beta^2\}(1 + |t|/\log^{100}x))} \left(\frac{\log y}{\log x} \right)^{\alpha^2 + \beta^2} \left(1 + \min\{\frac{\log y}{\log x}, \frac{1}{|t|\log x}\} \right)^{2\alpha \beta} . $$
\end{epres1}

\begin{proof}[Proof of Euler Product Result 1]
For concision in writing the proof, let us temporarily set $M = M(\alpha,\beta) := \max\{\alpha, \beta, \alpha^3 , \beta^3\}$.

Firstly we may rewrite
\begin{eqnarray}
&& \left|1 - \frac{f(p)}{p^{1/2+\sigma}}\right|^{-2\alpha} \left|1 - \frac{f(p)}{p^{1/2+\sigma + it}}\right|^{-2\beta} = \exp\{-2\alpha \Re\log(1 - \frac{f(p)}{p^{1/2+\sigma}}) -2\beta \Re\log(1 - \frac{f(p)}{p^{1/2+\sigma + it}})\} \nonumber \\
& = & \exp\{ \frac{2\alpha \Re f(p)}{p^{1/2 + \sigma}} + \frac{\alpha \Re f(p)^2}{p^{1+2\sigma}} +  \frac{2\beta \Re f(p)p^{-it}}{p^{1/2 + \sigma}} + \frac{\beta \Re f(p)^2 p^{-2it}}{p^{1+2\sigma}} + O(\frac{\max\{\alpha, \beta\}}{p^{3/2 + 3\sigma}}) \} . \nonumber
\end{eqnarray}
Next, if $y \geq p > x \geq 100\max\{\alpha^2 , \beta^2\}$ then every term in the exponential here has size at most $2\max\{\alpha, \beta\}/p^{1/2 + \sigma} = 2\max\{\alpha, \beta\} e^{-\sigma \log p}/ p^{1/2} \leq e/5$. Therefore we may apply the series expansion of the exponential function, finding the above is
\begin{eqnarray}
& = & 1 + \frac{2(\alpha \Re f(p) + \beta \Re f(p)p^{-it})}{p^{1/2 + \sigma}} + \frac{(\alpha \Re f(p)^2 + \beta \Re f(p)^2 p^{-2it})}{p^{1+2\sigma}} + \nonumber \\
&& + \frac{2(\alpha \Re f(p) + \beta \Re f(p)p^{-it})^2}{p^{1 + 2\sigma}} + O(\frac{M}{p^{3/2 + 3\sigma}}) . \nonumber
\end{eqnarray}
Now taking expectations, by symmetry we have $\E \Re f(p) = \E \Re f(p)^2 = 0$, similarly for $\E \Re f(p) p^{-it}$ and $\E \Re f(p)^2 p^{-2it}$. A simple trigonometric calculation also shows that $\E (\Re f(p))^2 = 1/2$, and similarly $\E \Re f(p) \Re f(p) p^{-it} = \cos(t\log p)/2$. So we get
\begin{eqnarray}
&& \E \left|1 - \frac{f(p)}{p^{1/2+\sigma}}\right|^{-2\alpha} \left|1 - \frac{f(p)}{p^{1/2+\sigma + it}}\right|^{-2\beta} = 1 + \frac{2 \E (\alpha \Re f(p) + \beta \Re f(p)p^{-it})^2}{p^{1 + 2\sigma}} + O(\frac{M}{p^{3/2 + 3\sigma}}) \nonumber \\
& = & 1 + \frac{2 (\alpha^2 \E (\Re f(p))^2 + 2\alpha\beta \E \Re f(p) \Re f(p) p^{-it} + \beta^2 \E (\Re f(p) p^{-it})^2)}{p^{1 + 2\sigma}} + O(\frac{M}{p^{3/2 + 3\sigma}}) \nonumber \\
& = & 1 + \frac{\alpha^2 + \beta^2 + 2\alpha\beta \cos(t\log p)}{p^{1 + 2\sigma}} + O(\frac{M}{p^{3/2 + 3\sigma}}) \nonumber \\
& = & \exp\{\frac{\alpha^2 + \beta^2 + 2\alpha\beta \cos(t\log p)}{p^{1 + 2\sigma}} + O(\frac{M}{p^{3/2 + 3\sigma}}) \} . \nonumber
\end{eqnarray}

Combining the above calculation with the independence of $f$ on distinct primes, and using that $p^{3/2 + \sigma} = e^{\sigma \log p} p^{3/2} \geq e^{-1} p^{3/2}$ for $p \leq y$, we deduce that the quantity $\E \prod_{x < p \leq y} \left|1 - \frac{f(p)}{p^{1/2+\sigma}}\right|^{-2\alpha} \left|1 - \frac{f(p)}{p^{1/2+\sigma + it}}\right|^{-2\beta}$ in the statement of the result is
\begin{eqnarray}
& = & \exp\{\sum_{x < p \leq y} (\frac{\alpha^2 + \beta^2 + 2\alpha\beta \cos(t\log p)}{p^{1 + 2\sigma}} + O(\frac{M}{p^{3/2 + 3\sigma}}) ) \} \nonumber \\
& = & \exp\{\sum_{x < p \leq y} \frac{\alpha^2 + \beta^2 + 2\alpha\beta \cos(t\log p)}{p^{1 + 2\sigma}} + O(\frac{\max\{\alpha. \beta, \alpha^3 , \beta^3\}}{\sqrt{x} \log x}) \} . \nonumber
\end{eqnarray}

To deduce the second part of Euler Product Result 1, we can use standard estimates from prime number theory. Indeed, the Chebychev and Mertens estimates for sums over primes imply that
\begin{eqnarray}
\sum_{x < p \leq y} \frac{\alpha^2 + \beta^2}{p^{1 + 2\sigma}} & = & (\alpha^2 + \beta^2) \sum_{x < p \leq y} \frac{1}{p} + (\alpha^2 + \beta^2) \sum_{x < p \leq y} \frac{e^{-2\sigma \log p} - 1}{p} \nonumber \\
& = & (\alpha^2 + \beta^2)\log(\frac{\log y}{\log x}) + O(\max\{\alpha^2 , \beta^2\}) , \nonumber
\end{eqnarray}
using that $e^{-2\sigma \log p} - 1 \ll |\sigma| \log p \ll \frac{\log p}{\log y}$ for $|\sigma| \leq 1/\log y$. We may remove the nuisance factor $p^{2\sigma}$ from the sum $\sum_{x < p \leq y} \frac{2\alpha\beta \cos(t\log p)}{p^{1 + 2\sigma}}$ with the same error term. Then using the Prime Number Theorem in the form $\pi(z) := \#\{p \leq z : p \; \text{prime}\} = \int_{2}^{z} \frac{du}{\log u} + O(\frac{z}{\log^{100}z})$, we have
\begin{eqnarray}
\sum_{x < p \leq y} \frac{\cos(t\log p)}{p} = \int_{x}^{y} \frac{\cos(t \log z)}{z} d\pi(z) & = & \int_{x}^{y} \frac{\cos(t \log z)}{z \log z} dz + O(\frac{1 + |t|}{\log^{100}x}) \nonumber \\
& = & \int_{\log x}^{\log y} \frac{\cos(tu)}{u} du + O(\frac{1+|t|}{\log^{100}x}) . \nonumber
\end{eqnarray}

Now if $|t|\log y \leq 1$, then the estimate $\cos(tu) = 1 + O((tu)^2)$ shows the integral is $\log\log y - \log\log x + O((t\log y)^2) = \log(\frac{\log y}{\log x}) + O(1)$. If instead we have $|t|\log x \leq 1$ but $|t|\log y > 1$, then we can evaluate the part of the integral with $u \leq 1/|t|$ using the estimate $\cos(tu) = 1 + O((tu)^2)$, and estimate the rest using integration by parts, yielding an overall estimate $\log(\frac{1}{|t|\log x}) + O(1)$. If $|t|\log x > 1$ then integration by parts shows the whole integral is $O(1)$. In any case, Euler Product Result 1 is proved.
\end{proof}

We will need a version of the above result for Rademacher random multiplicative functions. Unlike in the Steinhaus case, the distribution of $f(n)n^{-it}$ is not the same for all real $t$ in the Rademacher case, so our general statement must allow two different imaginary shifts in our two Euler product factors.

\begin{epres2}
If $f$ is a Rademacher random multiplicative function, then for any real $\alpha, \beta \geq 0$, any real $100(1 + \max\{\alpha^2 , \beta^2\}) \leq x \leq y$, and any real $\sigma \geq - 1/\log y$ and $t_1 , t_2$, we have
\begin{eqnarray}
&& \E \prod_{x < p \leq y} \left|1 + \frac{f(p)}{p^{1/2+\sigma + it_1}}\right|^{2\alpha} \left|1 + \frac{f(p)}{p^{1/2+\sigma + it_2}}\right|^{2\beta} \nonumber \\
& = & \exp\{\sum_{x < p \leq y} \frac{\alpha^2 + \beta^2 + (\alpha^2 - \alpha)\cos(2t_1 \log p) + (\beta^2 - \beta)\cos(2t_2 \log p)}{p^{1 + 2\sigma}} + \nonumber \\
&& + \sum_{x < p \leq y} \frac{2\alpha\beta(\cos((t_1 + t_2)\log p) + \cos((t_1 - t_2)\log p))}{p^{1 + 2\sigma}} + O(\frac{\max\{\alpha, \beta, \alpha^3 , \beta^3\}}{\sqrt{x} \log x}) \} . \nonumber
\end{eqnarray}

If we also have $\sigma \leq 1/\log y$, then the above is
\begin{eqnarray}
& = & e^{O(\max\{\alpha, \beta, \alpha^2 , \beta^2\}(1 + \frac{|t_1| + |t_2|}{\log^{100}x}))} \left(1 + \min\{\frac{\log y}{\log x}, \frac{|t_1|^{-1}}{\log x}\} \right)^{\alpha^2 - \alpha} \left(1 + \min\{\frac{\log y}{\log x}, \frac{|t_2|^{-1}}{\log x}\} \right)^{\beta^2 - \beta} \cdot \nonumber \\
&& \cdot \left(\frac{\log y}{\log x} \right)^{\alpha^2 + \beta^2} \left((1 + \min\{\frac{\log y}{\log x}, \frac{|t_1 + t_2|^{-1}}{\log x}\}) (1 + \min\{\frac{\log y}{\log x}, \frac{|t_1 - t_2|^{-1}}{\log x}\}) \right)^{2\alpha\beta} . \nonumber
\end{eqnarray}

As an upper bound, we may replace the error term $e^{O(\max\{\alpha, \beta, \alpha^2 , \beta^2\}(1 + \frac{|t_1| + |t_2|}{\log^{100}x}))}$ by $e^{O(\max\{\alpha, \beta, \alpha^2 , \beta^2\})} \min\{\frac{\log y}{\log x}, 1 + \frac{(|t_1| + |t_2|)^{1/100}}{\log x} \}^{|\alpha^2 - \alpha| + |\beta^2 - \beta| + 4\alpha\beta}$, and as a lower bound we may replace it by $e^{O(\max\{\alpha, \beta, \alpha^2 , \beta^2\})} \min\{\frac{\log y}{\log x}, 1 + \frac{(|t_1| + |t_2|)^{1/100}}{\log x} \}^{-(|\alpha^2 - \alpha| + |\beta^2 - \beta| + 4\alpha\beta)}$.
\end{epres2}

The estimation of the error terms here is rather crude, but will be sufficient as they only depend quite mildly on the $t_i$.

\begin{proof}[Proof of Euler Product Result 2]
The proof is a fairly straightforward adaptation of the proof of Euler Product Result 1. We again temporarily set $M = M(\alpha,\beta) := \max\{\alpha, \beta, \alpha^3 , \beta^3\}$. In the first place we have
\begin{eqnarray}
&& \left|1 + \frac{f(p)}{p^{1/2+\sigma+it_1}}\right|^{2\alpha} \left|1 + \frac{f(p)}{p^{1/2+\sigma + it_2}}\right|^{2\beta} \nonumber \\
& = & \exp\{2\alpha \Re\log(1 + \frac{f(p)}{p^{1/2+\sigma+it_1}}) +2\beta \Re\log(1 + \frac{f(p)}{p^{1/2+\sigma + it_2}})\} \nonumber \\
& = & \exp\{ \frac{2\alpha \Re f(p) p^{-it_1}}{p^{1/2 + \sigma}} - \frac{\alpha \Re f(p)^2 p^{-2it_1}}{p^{1+2\sigma}} +  \frac{2\beta \Re f(p)p^{-it_2}}{p^{1/2 + \sigma}} - \frac{\beta \Re f(p)^2 p^{-2it_2}}{p^{1+2\sigma}} + O(\frac{\max\{\alpha, \beta\}}{p^{3/2 + 3\sigma}}) \} \nonumber \\
& = & 1 + \frac{2(\alpha \Re f(p)p^{-it_1} + \beta \Re f(p)p^{-it_2})}{p^{1/2 + \sigma}} - \frac{(\alpha \Re f(p)^2 p^{-2it_1} + \beta \Re f(p)^2 p^{-2it_2})}{p^{1+2\sigma}} + \nonumber \\
&& + \frac{2(\alpha \Re f(p)p^{-it_1} + \beta \Re f(p)p^{-it_2})^2}{p^{1 + 2\sigma}} + O(\frac{M}{p^{3/2 + 3\sigma}}) . \nonumber
\end{eqnarray}
Furthermore, in the Rademacher case we have $f(p)^2 \equiv 1$, whilst still $\E \Re f(p) p^{-it} = \cos(t\log p) \E f(p) = 0$. So we get
\begin{eqnarray}
&& \E \left|1 + \frac{f(p)}{p^{1/2+\sigma+it_1}}\right|^{2\alpha} \left|1 + \frac{f(p)}{p^{1/2+\sigma + it_2}}\right|^{2\beta} \nonumber \\
& = & 1 - \frac{(\alpha\cos(2t_1 \log p) + \beta\cos(2t_2 \log p))}{p^{1+2\sigma}} + \frac{2 (\alpha\cos(t_1 \log p) + \beta\cos(t_2 \log p) )^2}{p^{1 + 2\sigma}} + O(\frac{M}{p^{3/2 + 3\sigma}}) , \nonumber
\end{eqnarray}
and using standard cosine identities this is all
\begin{eqnarray}
& = & 1 + \frac{\alpha^2 + \beta^2 + (\alpha^2 - \alpha)\cos(2t_1 \log p) + (\beta^2 - \beta)\cos(2t_2 \log p)}{p^{1 + 2\sigma}} + \nonumber \\
&& + \frac{2\alpha\beta(\cos((t_1 + t_2)\log p) + \cos((t_1 - t_2)\log p))}{p^{1 + 2\sigma}} + O(\frac{M}{p^{3/2 + 3\sigma}}) \nonumber \\
& = & \exp\{\frac{\alpha^2 + \beta^2 + (\alpha^2 - \alpha)\cos(2t_1 \log p) + (\beta^2 - \beta)\cos(2t_2 \log p)}{p^{1 + 2\sigma}} + \nonumber \\
&& + \frac{2\alpha\beta(\cos((t_1 + t_2)\log p) + \cos((t_1 - t_2)\log p))}{p^{1 + 2\sigma}} + O(\frac{M}{p^{3/2 + 3\sigma}}) \} . \nonumber
\end{eqnarray}

The first two conclusions of Euler Product Result 2 now follow exactly as in the proof of Euler Product Result 1.

For the final claimed inequalities, we note that the source of the unwanted error term $O(\max\{\alpha, \beta, \alpha^2 , \beta^2\} \frac{|t_1| + |t_2|}{\log^{100}x})$ in the exponent lies in our using the Prime Number Theorem to estimate the various sums $\sum_{x < p \leq y} \frac{(\alpha^2 - \alpha)\cos(2t_1 \log p)}{p^{1 + 2\sigma}}$, $\sum_{x < p \leq y} \frac{(\beta^2 - \beta)\cos(2t_2 \log p)}{p^{1 + 2\sigma}}$, $ \sum_{x < p \leq y} \frac{2\alpha\beta \cos((t_1 + t_2)\log p)}{p^{1 + 2\sigma}}$, $ \sum_{x < p \leq y} \frac{2\alpha\beta \cos((t_1 - t_2)\log p)}{p^{1 + 2\sigma}}$. Instead, if $|t_1| \geq \log^{100}x$ (which is the only case where it might produce a large error term) we can upper bound $\sum_{x < p \leq y} \frac{(\alpha^2 - \alpha)\cos(2t_1 \log p)}{p^{1 + 2\sigma}}$ by
$$ \sum_{x < p \leq \min\{e^{|t_1|^{1/100}},y\}} \frac{|\alpha^2 - \alpha|}{p^{1 + 2\sigma}} + \sum_{\min\{e^{|t_1|^{1/100}},y\} < p \leq y} \frac{(\alpha^2 - \alpha)\cos(2t_1 \log p)}{p^{1 + 2\sigma}} . $$
As in the proof of Euler Product Result 1, the second sum here is $\ll \max\{\alpha,\alpha^2\}$ (we can use the Prime Number Theorem to estimate it, since the lower end point is now sufficiently large that we don't pick up a big error term), and the first sum is
$$ |\alpha^2 - \alpha|(\sum_{x < p \leq \min\{e^{|t_1|^{1/100}},y\}} \frac{1}{p} + O(1)) = |\alpha^2 - \alpha| (\min\{\log(\frac{\log y}{\log x}), \log(\frac{|t_1|^{1/100}}{\log x})\} + O(1)) . $$
We can handle the other sums similarly when $t_2, t_1+t_2, t_1-t_2$ are large. In the worst case, as an upper bound this will produce an extra multiplicative factor
$$ \exp\{(|\alpha^2 - \alpha| + |\beta^2 - \beta| + 4\alpha\beta) \min\{\log(\frac{\log y}{\log x}), \log(1 + \frac{(|t_1|+|t_2|)^{1/100}}{\log x})\} \} . $$
An exactly similar argument gives a lower bound with $(|\alpha^2 - \alpha| + |\beta^2 - \beta| + 4\alpha\beta)$ replaced by $-(|\alpha^2 - \alpha| + |\beta^2 - \beta| + 4\alpha\beta)$.
\end{proof}

\subsection{Probabilistic preparations}
Next we record some moment estimates, mostly fairly simple yet interesting, that will be input to our arguments in various places.

\begin{probres1}[Rough hypercontractive inequalities]
For any real $q \geq 1$, the following is true.

If $f(n)$ is a Steinhaus random multiplicative function, then for any sequence of complex numbers $(a_{n})_{n \leq N}$ we have
$$ \E\left|\sum_{n \leq N} a_n f(n) \right|^{2q} \leq \left(\sum_{n \leq N} |a_n|^2 d_{\lceil q \rceil}(n) \right)^{q} , $$
where $d_{k}(\cdot)$ denotes the $k$-fold divisor function (i.e. the number of $k$-tuples of natural numbers whose product is $\cdot$, or equivalently the Dirichlet series coefficient of $\zeta(s)^k$), and $\lceil q \rceil$ denotes the ceiling of $q$.

If $f(n)$ is a Rademacher random multiplicative function, then for any sequence of complex numbers $(a_{n})_{n \leq N}$ we have
$$ \E\left|\sum_{n \leq N} a_n f(n) \right|^{2q} \leq \left(\sum_{n \leq N} |a_n|^2 d_{2\lceil q \rceil - 1}(n) \right)^{q} . $$
\end{probres1}

\begin{proof}[Proof of Probability Result 1]
By H\"{o}lder's inequality, it suffices to treat the case where $q$ is a natural number.

For Steinhaus $f(n)$, expanding the $2q$-th power and taking expectations we get
$$ \E\left|\sum_{n \leq N} a_n f(n) \right|^{2q} = \sum_{n_1, ..., n_q \leq N} a_{n_1} ... a_{n_{q}} \sum_{m_1, ..., m_q \leq N} \overline{a_{m_1} ... a_{m_{q}}} \textbf{1}_{\prod_{i=1}^{q} n_i = \prod_{i=1}^{q} m_i} , $$
where $\textbf{1}$ denotes the indicator function. Using the upper bound $|a_{n_1} ... a_{n_{q}} \overline{a_{m_1} ... a_{m_{q}}}| \leq (1/2)(|a_{n_1} ... a_{n_{q}}|^2 + |a_{m_1} ... a_{m_{q}}|^2)$, together with the symmetry of the $n_i$ and the $m_i$, we deduce that $\E\left|\sum_{n \leq N} a_n f(n) \right|^{2q}$ is
$$ \leq \sum_{n_1, ..., n_q \leq N} |a_{n_1} ... a_{n_{q}}|^2 \sum_{m_1, ..., m_q \leq N} \textbf{1}_{\prod_{i=1}^{q} n_i = \prod_{i=1}^{q} m_i} \leq \sum_{n_1, ..., n_q \leq N} |a_{n_1} ... a_{n_{q}}|^2 d_{q}\left(\prod_{i=1}^{q} n_i \right) . $$
Finally, since the function $d_{q}(\cdot)$ is sub-multiplicative we find the above is
$$ \leq \sum_{n_1, ..., n_q \leq N} |a_{n_1} ... a_{n_{q}}|^2 d_{q}(n_1) ... d_{q}(n_q) = \left(\sum_{n \leq N} |a_n|^2 d_{q}(n) \right)^{q} . $$

In the Rademacher case, one needs a bit more involved argument. We refer the reader to Lemma 2 of Hal\'{a}sz~\cite{halasz}, where this result is proved by induction on the exponent $2q$. We may remark that, since Rademacher $f(n)$ is only supported on squarefree $n$, we may assume that $a_n$ is only non-zero for squarefree $n$, and then $d_{2\lceil q \rceil - 1}(n) = (2\lceil q \rceil - 1)^{\Omega(n)}$ where $\Omega(n)$ is the number of prime factors of $n$. The ultimate source of the factors $d_{2\lceil q \rceil - 1}(n)$ is that, when one expands the expectation in the inductive proof, the only surviving terms are those where the product $n_1 ... n_{2q}$ is a perfect square, so all the prime factors of $n_{2q}$ must be repeated somewhere amongst the other terms $n_1, ..., n_{2q-1}$.
\end{proof}

We describe the inequalities in Probability Result 1 as ``rough hypercontractive inequalities'' because (if we take $2q$-th roots of both sides) they upper bound an $L^{2q}$-norm by a weighted $L^{2}$ norm without any other terms, but the weights $d_{\lceil q \rceil}(n), d_{2\lceil q \rceil - 1}(n)$ will not generally be the sharpest possible unless $q$ is an integer. One can prove more precise results for non-integer $q$ using more subtle interpolation techniques, see section 2 of Bondarenko, Brevig, Saksman, Seip and Zhao~\cite{bbsszpseudo} for the Steinhaus case, and Chapitre III of Bonami~\cite{bonami} for the Rademacher case (expressed in rather different notation). However, for our applications the extra precision in these inequalities will not be needed.

\begin{probres2}
Let $(\epsilon_{n})_{n \leq N}$ be a sequence of independent random variables, each satisfying $\E \epsilon_n = 0$ and $\E |\epsilon_n|^2 = 1$, and let $(a_{n})_{n \leq N}$ be a sequence of complex numbers. Then for any real $q \geq 1$, we have
$$ \E\left| \sum_{n \leq N} a_n \epsilon_n \right|^{2q} \geq \left( \sum_{n \leq N} |a_n|^2 \right)^{q} . $$
\end{probres2}

\begin{proof}[Proof of Probability Result 2]
Since we assume that $q \geq 1$, simply applying H\"{o}lder's inequality we get
$$ \left( \sum_{n \leq N} |a_n|^2 \right)^{q} = \left( \E\left| \sum_{n \leq N} a_n \epsilon_n \right|^{2} \right)^{q} \leq \E\left| \sum_{n \leq N} a_n \epsilon_n \right|^{2q} . $$
\end{proof}

If the $\epsilon_n$ are Rademacher or Steinhaus random variables\footnote{We emphasise that here we are referring to ordinary Rademacher or Steinhaus random variables, not random {\em multiplicative} functions.}, then Khintchine's inequality (see e.g. Lemma 3.8.1 of Gut~\cite{gut}) in fact implies that $\E\left| \sum_{n \leq N} a_n \epsilon_n \right|^{2q} \asymp_{q} \left( \sum_{n \leq N} |a_n|^2 \right)^{q}$ for all real $q \geq 0$. For our purposes here we will only require the simple lower bound in Probability Result 2, but it is useful to keep Khintchine's inequality in mind since it means that when we apply the lower bound, we are doing something sharp.

The final result we shall record is more sophisticated, and requires some terminology before we can state it. Suppose that $(\Omega, \mathcal{F}, \p)$ is a probability space, and $(\mathcal{F}_{n})_{n \geq 0}$ is a {\em filtration} on $\mathcal{F}$, in other words a sequence of sub-$\sigma$-algebras satisfying $\mathcal{F}_{0} \subseteq \mathcal{F}_{1} \subseteq ... \subseteq \mathcal{F}$. We say a sequence of random variables $(X_n)_{n \geq 0}$ on $(\Omega, \mathcal{F}, \p)$ is a {\em submartingale} (relative to $(\mathcal{F}_{n})_{n \geq 0}$ and $\p$) if it satisfies:
\begin{enumerate}
\item (adapted) $X_n$ is measurable with respect to $\mathcal{F}_{n}$, for all $n \geq 0$;

\item (integrable) $\E|X_n|$ is finite, for all $n \geq 0$;

\item (non-decreasing on average) for all $n \geq 1$, the conditional expectation $\E(X_n |\mathcal{F}_{n-1}) \geq X_{n-1}$ almost surely.
\end{enumerate}

Condition (iii) says that a submartingale is non-decreasing on average, in quite a strong sense: for any given value of $X_{n-1}$ (or, informally speaking, any other ``information'' from the sigma algebra $\mathcal{F}_{n-1}$), the conditional expectation of $X_n$ will be at least as large. One can apply this property to partition the sample space $\Omega$ in useful ways, and prove that the moments of the random variables comprising a submartingale satisfy the following useful bound. We will use this as an ingredient in proving our $2q$-th moment upper bounds when $q$ is close to 1.
\begin{probres3}[Doob's $L^p$ maximal inequality, see Theorem 9.4 of Gut~\cite{gut}]
Let $(X_n)_{n \geq 0}$ be a non-negative submartingale (on some probability space and with respect to some filtration). Then for any $p > 1$, we have
$$ \E(\max_{0 \leq k \leq n} X_k)^p \leq (\frac{p}{p-1})^p \E X_{n}^p . $$
\end{probres3}

\subsection{Some miscellaneous lemmas}
As in the first paper~\cite{harperrmflow} in this sequence, we will need the following version of Parseval's identity for Dirichlet series to help with relating $\E|\sum_{n \leq x} f(n)|^{2q}$ to an Euler product average.

\begin{harman1}[See (5.26) in sec. 5.1 of Montgomery and Vaughan~\cite{mv}]
Let $(a_n)_{n=1}^{\infty}$ be any sequence of complex numbers, and let $A(s) := \sum_{n=1}^{\infty} \frac{a_n}{n^s}$ denote the corresponding Dirichlet series, and $\sigma_c$ denote its abscissa of convergence. Then for any $\sigma > \max\{0,\sigma_c \}$, we have
$$ \int_{0}^{\infty} \frac{|\sum_{n \leq x} a_n |^2}{x^{1 + 2\sigma}} dx = \frac{1}{2\pi} \int_{-\infty}^{\infty} \left|\frac{A(\sigma + it)}{\sigma + it}\right|^2 dt . $$
\end{harman1}

We will use the following estimate to handle sums of divisor-type functions that appear in our calculations.

\begin{numth1}[See Lemma 2.1 of Lau, Tenenbaum and Wu~\cite{tenenbaum}]
Let $0 < \delta < 1$, let $m \geq 1$, and suppose that $\max\{3, 2m\} \leq y \leq z \leq y^{10}$ and that $1 < u \leq v(1 - y^{-\delta})$. As usual, let $\Omega(d)$ denote the total number of prime factors of $d$ (counted with multiplicity). Then
$$ \sum_{\substack{u \leq d \leq v, \\ p | d \Rightarrow y \leq p \leq z}} m^{\Omega(d)} \ll_{\delta} \frac{(v-u) m}{\log y} \prod_{y \leq p \leq z} \left(1 - \frac{m}{p}\right)^{-1} . $$
\end{numth1}

This is a slight generalisation of Lemma 2.1 of Lau, Tenenbaum and Wu~\cite{tenenbaum} (see also Lemma 3 of Hal\'{a}sz~\cite{halasz}). See section 2.1 of Harper~\cite{harperrmflow} for the full (short) proof.

\section{Easier cases of the theorems}\label{seceasycases}
As remarked in the Introduction, since we allow a multiplicative error term $e^{O(q^2)}$ in our Theorems, it turns out that proving our claimed upper bounds when $\log\log x \leq q \leq \frac{c\log x}{\log\log x}$ is somewhat straightforward. We present these arguments in this section. Some of the techniques involved, including the use of Rankin's trick with an exponent roughly like $1 + q/\log x$, and a special treatment of prime factors that are $\ll q^2$, will recur later when we develop our main arguments.

\subsection{The upper bound in the Steinhaus case, for very large $q$}\label{secsteineasy}
For $q \geq \log\log x$ we have $\log^{(q-1)^2}x = \log^{q^2 + O(q)}x = e^{O(q^2)} \log^{q^2}x$. Thus to establish the upper bound part of Theorem 1 for $\log\log x \leq q \leq \frac{c\log x}{\log\log x}$, it will suffice to show that
$$ \vert\vert \sum_{n \leq x} f(n) \vert\vert_{2q} \leq e^{-(q/2)\log q - (q/2)\log\log(2q) + O(q)} \sqrt{x} \log^{q/2}x , $$
where as usual we write $\vert\vert \cdot \vert\vert_{r} := (\E|\cdot|^{r})^{1/r}$.

To prove this, we first apply Minkowski's inequality to obtain that
$$ \vert\vert \sum_{n \leq x} f(n) \vert\vert_{2q} = \vert\vert \sum_{\substack{m \leq x, \\ m \; \text{is} \; q^2 \; \text{smooth}}} f(m) \sum_{\substack{n \leq x/m, \\ p \mid n \Rightarrow p > q^2}} f(n) \vert\vert_{2q} \leq \sum_{\substack{m \leq x, \\ m \; \text{is} \; q^2 \; \text{smooth}}}  \vert\vert \sum_{\substack{n \leq x/m, \\ p \mid n \Rightarrow p > q^2}} f(n) \vert\vert_{2q} . $$
Recall here that a number is said to be $q^2$-smooth if all of its prime factors are $\leq q^2$. Using the first part of Probability Result 1, and then using Rankin's trick of upper bounding $\textbf{1}_{n \leq x/m}$ by $(\frac{x}{nm})^{1+q/\log x}$ (and recalling that the divisor function $d_{\lceil q \rceil}(n)$ is the Dirichlet series coefficient of $\zeta(s)^{\lceil q \rceil} = \sum_{n=1}^{\infty} \frac{d_{\lceil q \rceil}(n)}{n^s} = \prod_{p}(1 - \frac{1}{p^s})^{-\lceil q \rceil}$), we get
\begin{eqnarray}
\vert\vert \sum_{\substack{n \leq x/m, \\ p \mid n \Rightarrow p > q^2}} f(n) \vert\vert_{2q} \leq \Biggl (\sum_{\substack{n \leq x/m, \\ p \mid n \Rightarrow p > q^2}} d_{\lceil q \rceil}(n) \Biggr)^{1/2} & \leq & \Biggl ( (\frac{x}{m})^{1 + q/\log x} \sum_{\substack{n \leq x/m, \\ p \mid n \Rightarrow p > q^2}} \frac{d_{\lceil q \rceil}(n)}{n^{1+q/\log x}} \Biggr)^{1/2} \nonumber \\
& \leq & \sqrt{\frac{x}{m}} e^{O(q)} \prod_{p > q^2}(1 - \frac{1}{p^{1 + q/\log x}})^{- \lceil q \rceil/2} . \nonumber
\end{eqnarray}

Finally, the product over primes here is $\zeta(1 + \frac{q}{\log x})^{\lceil q \rceil/2} \prod_{p \leq q^2}(1 - \frac{1}{p^{1 + q/\log x}})^{\lceil q \rceil/2}$, which is $= e^{O(q)} (\frac{\log x}{q})^{\lceil q \rceil/2} e^{ - \sum_{p \leq q^2} \frac{\lceil q \rceil/2}{p^{1 + q/\log x}}} = e^{O(q)} (\frac{\log x}{q \log q})^{q/2}$ on our range $\log\log x \leq q \leq \frac{c\log x}{\log\log x}$. And when we sum over $m$ we have $\sum_{\substack{m \leq x, \\ m \; \text{is} \; q^2 \; \text{smooth}}} \frac{1}{\sqrt{m}} \leq e^{\sum_{p \leq q^2} O(1/\sqrt{p})} \leq e^{O(q)} $, so putting everything together we get an acceptable upper bound for $\vert\vert \sum_{n \leq x} f(n) \vert\vert_{2q}$.
\qed

\subsection{The upper bound in the Rademacher case, for very large $q$}
Similarly as in the Steinhaus case, to prove the upper bound part of Theorem 2 for $\log\log x \leq q \leq \frac{c\log x}{\log\log x}$ it will suffice to show that, for Rademacher random multiplicative $f(n)$, we have
$$ \vert\vert \sum_{n \leq x} f(n) \vert\vert_{2q} \leq e^{-q\log q - q\log\log(2q) + O(q)} \sqrt{x} \log^{q}x . $$

Using Minkowski's inequality, the second part of Probability Result 1, and then Rankin's trick, we get
\begin{eqnarray}
\vert\vert \sum_{n \leq x} f(n) \vert\vert_{2q} & \leq & \sum_{\substack{m \leq x, \\ m \; \text{is} \; q^2 \; \text{smooth}}}  \vert\vert \sum_{\substack{n \leq x/m, \\ p \mid n \Rightarrow p > q^2}} f(n) \vert\vert_{2q} \leq \sum_{\substack{m \leq x, \\ m \; \text{is} \; q^2 \; \text{smooth}}} \Biggl (\sum_{\substack{n \leq x/m, \\ p \mid n \Rightarrow p > q^2}} d_{2\lceil q \rceil - 1}(n) \Biggr)^{1/2} \nonumber \\
& \leq & \sum_{\substack{m \leq x, \\ m \; \text{is} \; q^2 \; \text{smooth}}} \sqrt{\frac{x}{m}} e^{O(q)} \prod_{p > q^2}(1 - \frac{1}{p^{1 + q/\log x}})^{- (2\lceil q \rceil - 1)/2} . \nonumber
\end{eqnarray}

We can estimate the product over primes as in the Steinhaus case, finding it is $= e^{O(q)} (\frac{\log x}{q})^{(2\lceil q \rceil - 1)/2} e^{ - \sum_{p \leq q^2} \frac{(2\lceil q \rceil - 1)/2}{p^{1 + q/\log x}}} = e^{O(q)} (\frac{\log x}{q \log q})^{q}$ on our range $\log\log x \leq q \leq \frac{c\log x}{\log\log x}$. And when we sum over $m$ we again have $\sum_{\substack{m \leq x, \\ m \; \text{is} \; q^2 \; \text{smooth}}} \frac{1}{\sqrt{m}} \leq e^{\sum_{p \leq q^2} O(1/\sqrt{p})} \leq e^{O(q)} $, so putting everything together we get an acceptable upper bound for $\vert\vert \sum_{n \leq x} f(n) \vert\vert_{2q}$.
\qed

\section{The reduction to Euler products}\label{seceulerprod}
In this section we shall prove four Propositions that make precise the assertion in \eqref{introprod}, that $\E|\sum_{n \leq x} f(n)|^{2q}$ may be bounded by studying integrals of Euler products.

\subsection{Upper bounds: statement of the propositions}
We will need a little notation, which is exactly the same as in the author's previous paper~\cite{harperrmflow} dealing with low moments. Given a random multiplicative function $f(n)$ (either Steinhaus or Rademacher, depending on the context), and an integer $0 \leq k \leq \log\log x$, let $F_k$ denote the partial Euler product of $f(n)$ over $x^{e^{-(k+1)}}$-smooth numbers. Thus for all complex $s$ with $\Re(s) > 0$, we have
$$ F_{k}(s) = \prod_{p \leq x^{e^{-(k+1)}}} \left(1 - \frac{f(p)}{p^s}\right)^{-1} = \sum_{\substack{n=1, \\ n \; \text{is} \; x^{e^{-(k+1)}} \; \text{smooth}}}^{\infty} \frac{f(n)}{n^s} $$
in the Steinhaus case, and
$$ F_{k}(s) = \prod_{p \leq x^{e^{-(k+1)}}} \left(1 + \frac{f(p)}{p^s}\right) = \sum_{\substack{n=1, \\ n \; \text{is} \; x^{e^{-(k+1)}} \; \text{smooth}}}^{\infty} \frac{f(n)}{n^s} $$
in the Rademacher case (the product taking a different form because $f(n)$ is only supported on squarefree numbers in that case).

\begin{prop1}\label{propstupper}
Let $f(n)$ be a Steinhaus random multiplicative function, let $x$ be large, and set $\mathcal{L} := \lfloor (\log\log x)/10 \rfloor$. Uniformly for all $1 \leq q \leq \log^{0.05}x$, we have
$$ \vert\vert \sum_{n \leq x} f(n) \vert\vert_{2q} \leq \sqrt{\frac{x}{\log x}} e^{O(q)} \sum_{0 \leq k \leq \mathcal{L}} \vert\vert \int_{-1/2}^{1/2} |F_k(1/2 + \frac{q-k}{\log x} + it)|^2 dt \vert\vert_{q}^{1/2} + e^{O(q)} \sqrt{\frac{x}{\log x}} . $$
\end{prop1}

In the low moments case, Proposition 1 of Harper~\cite{harperrmflow} gives an analogous upper bound for all $2/3 \leq q \leq 1$, but with the quantity $\mathcal{L}$ replaced by the smaller quantity $\mathcal{K} = \lfloor \log\log\log x \rfloor$, and the shift $\frac{q-k}{\log x}$ in the Euler product replaced by $\frac{-k}{\log x}$.

The additional shift by $\frac{q}{\log x}$ here corresponds to applying Rankin's trick with exponent $1 + q/\log x$ in our treatment of very large $q$ in section \ref{seceasycases}. We can introduce this at the acceptable cost of a prefactor $e^{O(q)}$ in the proposition, and it means that when we analyse the Euler product we can restrict attention to numbers that are $x^{1/q}$-smooth, which is crucial to obtaining the desired factor $e^{-q^2 \log q}$ in Theorem 1. The significant contribution from very smooth numbers, when $q$ becomes large, also explains why we must let $k$ run over a wider range than in the low moments case to obtain acceptable bounds. Finally, we remark that the range $1 \leq q \leq \log^{0.05}x$ allowed in Proposition 1 is somewhat artificial, but more than sufficient since we already proved the Theorem 1 upper bound for all $\log\log x \leq q \leq \frac{c\log x}{\log\log x}$ in section \ref{seceasycases}. It could be increased somewhat, but it seems hard to obtain an upper bound of a similar shape to Proposition 1 on the full range $1 \leq q \leq \frac{c\log x}{\log\log x}$, since for very large $q$ the significant contribution from very smooth numbers changes the behaviour in parts of the proof.

\begin{prop2}\label{propradupper}
Let $f(n)$ be a Rademacher random multiplicative function, let $x$ be large, and set $\mathcal{L} := \lfloor (\log\log x)/10 \rfloor$. Uniformly for all $1 \leq q \leq \log^{0.05}x$, we have
\begin{eqnarray}
\vert\vert \sum_{n \leq x} f(n) \vert\vert_{2q} & \leq & \sqrt{\frac{x}{\log x}} e^{O(q)} \sum_{0 \leq k \leq \mathcal{L}} \max_{N \in \Z} \frac{1}{(|N|+1)^{1/8}} \vert\vert \int_{N-1/2}^{N+1/2} |F_k(1/2 + \frac{q-k}{\log x} + it)|^2 dt \vert\vert_{q}^{1/2} \nonumber \\
&& + e^{O(q)} \sqrt{\frac{x}{\log x}} . \nonumber
\end{eqnarray}
\end{prop2}

One has to deal with translates by $N$ in the Rademacher case because, unlike in the Steinhaus case, the distribution of $(f(n)n^{it})$ is {\em not} the same (for $t \neq 0$) as the distribution of $(f(n))$ for Rademacher random multiplicative $f(n)$. However, as in the low moments argument in ~\cite{harperrmflow}, the main contribution will come from small $N$.

\subsection{Lower bounds: statement of the propositions}
For our work on lower bounds, we again connect the size of $\vert\vert \sum_{n \leq x} f(n) \vert\vert_{2q}$ with a certain integral average, and thence with random Euler products. Let $F$ denote the partial Euler product of $f(n)$, either Steinhaus or Rademacher, over $x$-smooth numbers. (Thus $F = F_{-1}$, if we slightly abuse our earlier notation).

\begin{prop3}\label{propstlower}
If $f(n)$ is a Steinhaus random multiplicative function, and $x$ is large, then uniformly for all $q \geq 1$ we have
$$ \vert\vert \sum_{n \leq x} f(n) \vert\vert_{2q} \gg \sqrt{\frac{x}{\log x}} \vert\vert \int_{1}^{x^{1/4}} \left|\sum_{m \leq z} f(m) \right|^2 \frac{dz}{z^{2}} \vert\vert_{q}^{1/2} . $$

In particular, for any large quantity $V \leq (\log x)/q$ we have that $\vert\vert \sum_{n \leq x} f(n) \vert\vert_{2q}$ is
$$ \gg \sqrt{\frac{x}{\log x}} \Biggl( \vert\vert \int_{-1/2}^{1/2} |F(1/2 + \frac{4Vq}{\log x} + it)|^2 dt \vert\vert_{q}^{1/2} - \frac{C}{e^{Vq/2}} \vert\vert \int_{-1/2}^{1/2} |F(1/2 + \frac{2Vq}{\log x} + it)|^2 dt \vert\vert_{q}^{1/2} \Biggr) , $$
where $C > 0$ is an absolute constant.
\end{prop3}

Notice that we don't need to impose any upper bound on $q$ here (although, for the second statement, there is an implicit upper bound $q \ll \log x$ in order that we can choose large $V \leq (\log x)/q$). This means we can use Proposition 3 to prove the lower bound in Theorem 1 on the full range of $q$ there.

\begin{prop4}\label{propradlower}
If $f(n)$ is a Rademacher random multiplicative function, the first bound in Proposition 3 continues to hold, and the second bound may be replaced by the statement that
\begin{eqnarray}
\vert\vert \sum_{n \leq x} f(n) \vert\vert_{2q} & \gg & \sqrt{\frac{x}{\log x}} \Biggl( \vert\vert \int_{-1/2}^{1/2} |F(1/2 + \frac{4Vq}{\log x} + it)|^2 dt \vert\vert_{q}^{1/2} - \nonumber \\
&& - \frac{C}{e^{Vq/2}} \max_{N \in \Z} \frac{1}{(|N|+1)^{1/8}} \vert\vert \int_{N-1/2}^{N+1/2} |F(1/2 + \frac{2Vq}{\log x} + it)|^2 dt \vert\vert_{q}^{1/2} \Biggr) . \nonumber
\end{eqnarray}
\end{prop4}

These results are again of the same general shape as the corresponding Propositions 3 and 4 of Harper~\cite{harperrmflow} from the low moments case. In fact, the propositions here are a little simpler as they don't involve an additional subtracted error term $-C \sqrt{\frac{x}{\log x}}$. This is accomplished by some reorganisation of the proof, and shrinking the range of integration over $z$ to $[1,x^{1/4}]$ rather than $[1,\sqrt{x}]$ from the low moments case, which makes no difference when applying the results. The other difference, similarly as in section \ref{seceasycases} and in our discussion of upper bounds, is that here we introduce shifts of the shape $\frac{4Vq}{\log x}$ in our Euler products, as opposed to $\frac{4V}{\log x}$ in the low moments analogues.

\subsection{Proof of Propositions 1 and 2}
We begin with Proposition 1. Let $P(n)$ denote the largest prime factor of $n$, and recall that a number $n$ is said to be $y$-smooth if $P(n) \leq y$. Recall also that the divisor function $d_{\lceil q \rceil}(n)$ is the Dirichlet series coefficient of $\zeta(s)^{\lceil q \rceil} = \prod_{p}(1 - \frac{1}{p^s})^{-\lceil q \rceil}$. By Minkowski's inequality, we have
$$ \vert\vert \sum_{n \leq x} f(n) \vert\vert_{2q} \leq \sum_{0 \leq k \leq \mathcal{L}} \vert\vert \sum_{\substack{n \leq x, \\ x^{e^{-(k+1)}} < P(n) \leq x^{e^{-k}}}} f(n) \vert\vert_{2q} + \vert\vert \sum_{\substack{n \leq x, \\ P(n) \leq x^{e^{-(\mathcal{L}+1)}} }} f(n) \vert\vert_{2q} . $$
Furthermore, the first part of Probability Result 1, followed by Rankin's trick with exponent $1- \frac{1}{\log^{0.9}x}$ (bounding $\textbf{1}_{n \leq x}$ by $(\frac{x}{n})^{1-1/\log^{0.9}x} = xe^{-\log^{0.1}x} \frac{1}{n^{1-1/\log^{0.9}x}}$), imply that
$$ \vert\vert \sum_{\substack{n \leq x, \\ P(n) \leq x^{e^{-(\mathcal{L}+1)}} }} f(n) \vert\vert_{2q} \leq \sqrt{\sum_{\substack{n \leq x, \\ P(n) \leq x^{e^{-(\mathcal{L}+1)}} }} d_{\lceil q \rceil}(n) } \leq \sqrt{x e^{-\log^{0.1}x} \sum_{\substack{n \leq x, \\ P(n) \leq x^{e^{-(\mathcal{L}+1)}} }} \frac{d_{\lceil q \rceil}(n)}{n^{1 - 1/\log^{0.9}x}} } . $$
Here the sum over $n$ is $\leq \prod_{p \leq x^{e^{-(\mathcal{L}+1)}}} (1 - \frac{1}{p^{1 - 1/\log^{0.9}x}})^{-\lceil q \rceil}$, and recalling that $\mathcal{L} := \lfloor (\log\log x)/10 \rfloor$ this is $\leq \prod_{p \leq e^{\log^{0.9}x}} (1 - \frac{1}{p^{1 - 1/\log^{0.9}x}})^{-\lceil q \rceil}$, which is $= e^{O(q)} \prod_{p \leq e^{\log^{0.9}x}} (1 - \frac{1}{p})^{-\lceil q \rceil} = \log^{O(q)}x$ by standard Chebychev and Mertens estimates for sums over primes. Since we assume in Proposition 1 that $q \leq \log^{0.05}x$, this whole contribution is $\ll \sqrt{x} e^{-c\log^{0.1}x}$, which is more than acceptable.

Next, if we let $\E^{(k)}$ denote expectation conditional on $(f(p))_{p \leq x^{e^{-(k+1)}}}$, then the first part of Probability Result 1 (applied, after conditioning on $(f(p))_{p \leq x^{e^{-(k+1)}}}$, with $a_m = \textbf{1}_{p|m \Rightarrow x^{e^{-(k+1)}} < p \leq x^{e^{-k}}} \cdot \sum_{\substack{n \leq x/m, \\ n \; \text{is} \; x^{e^{-(k+1)}} \text{-smooth}}} f(n)$) implies $\sum_{0 \leq k \leq \mathcal{L}} \vert\vert \sum_{\substack{n \leq x, \\ x^{e^{-(k+1)}} < P(n) \leq x^{e^{-k}}}} f(n) \vert\vert_{2q}$ is
\begin{eqnarray}
& = & \sum_{0 \leq k \leq \mathcal{L}} \vert\vert \sum_{\substack{1 < m \leq x , \\ p|m \Rightarrow x^{e^{-(k+1)}} < p \leq x^{e^{-k}}}} f(m) \sum_{\substack{n \leq x/m, \\ n \; \text{is} \; x^{e^{-(k+1)}} \text{-smooth}}} f(n)  \vert\vert_{2q} . \nonumber \\
& = & \sum_{0 \leq k \leq \mathcal{L}} \Biggl( \E \E^{(k)}\Biggl|\sum_{\substack{1 < m \leq x , \\ p|m \Rightarrow x^{e^{-(k+1)}} < p \leq x^{e^{-k}}}} f(m) \sum_{\substack{n \leq x/m, \\ n \; \text{is} \; x^{e^{-(k+1)}} \text{-smooth}}} f(n) \Biggr|^{2q} \Biggr)^{1/2q} \nonumber \\
& \leq & \sum_{0 \leq k \leq \mathcal{L}} \Biggl( \E( \sum_{\substack{1 < m \leq x , \\ p|m \Rightarrow x^{e^{-(k+1)}} < p \leq x^{e^{-k}}}} d_{\lceil q \rceil}(m) \Biggl|\sum_{\substack{n \leq x/m, \\ n \; \text{is} \; x^{e^{-(k+1)}} \text{-smooth}}} f(n) \Biggr|^2)^{q} \Biggr)^{1/2q} . \nonumber
\end{eqnarray}

To proceed further, we want to replace $\Biggl|\sum_{\substack{n \leq x/m, \\ n \; \text{is} \; x^{e^{-(k+1)}} \text{-smooth}}} f(n) \Biggr|^2$ in the above by a smoothed version. Set $X = e^{\sqrt{\log x}}$, say, and note that (uniformly for any $1 \leq q \leq \log^{0.05}x$) the above is
\begin{eqnarray}\label{uppersmoothingdisplay}
& = & \sum_{0 \leq k \leq \mathcal{L}} \vert\vert \sum_{\substack{1 < m \leq x , \\ p|m \Rightarrow x^{e^{-(k+1)}} < p \leq x^{e^{-k}}}} d_{\lceil q \rceil}(m) \Biggl|\sum_{\substack{n \leq x/m, \\ n \; \text{is} \; x^{e^{-(k+1)}} \text{-smooth}}} f(n) \Biggr|^2 \vert\vert_{q}^{1/2} \nonumber \\
& \ll & \sum_{0 \leq k \leq \mathcal{L}} \vert\vert \sum_{\substack{1 < m \leq x , \\ p|m \Rightarrow x^{e^{-(k+1)}} < p \leq x^{e^{-k}}}} d_{\lceil q \rceil}(m) \frac{X}{m} \int_{m}^{m(1 + \frac{1}{X})} \Biggl| \sum_{\substack{n \leq x/t, \\ x^{e^{-(k+1)}} \text{-smooth}}} f(n) \Biggr|^2 dt \vert\vert_{q}^{1/2} \nonumber \\
&& + \sum_{0 \leq k \leq \mathcal{L}} \vert\vert \sum_{\substack{1 < m \leq x , \\ p|m \Rightarrow x^{e^{-(k+1)}} < p \leq x^{e^{-k}}}} d_{\lceil q \rceil}(m) \frac{X}{m} \int_{m}^{m(1 + \frac{1}{X})} \Biggl| \sum_{\substack{x/t < n \leq x/m, \\ x^{e^{-(k+1)}} \text{-smooth}}} f(n) \Biggr|^2 dt \vert\vert_{q}^{1/2} .
\end{eqnarray}

We next want to show that the second term in \eqref{uppersmoothingdisplay} may be discarded as an error term. Using Minkowski's inequality again, followed by H\"{o}lder's inequality with exponent $q$ applied to the normalised integral $\frac{X}{m} \int_{m}^{m(1 + \frac{1}{X})} dt$, this second term is
\begin{eqnarray}
& \leq & \sum_{0 \leq k \leq \mathcal{L}} \sqrt{\sum_{\substack{1 < m \leq x , \\ p|m \Rightarrow x^{e^{-(k+1)}} < p \leq x^{e^{-k}}}} d_{\lceil q \rceil}(m) \vert\vert \frac{X}{m} \int_{m}^{m(1 + \frac{1}{X})} \Biggl| \sum_{\substack{x/t < n \leq x/m, \\ x^{e^{-(k+1)}} \text{-smooth}}} f(n) \Biggr|^2 dt \vert\vert_{q} } \nonumber \\
& \leq & \sum_{0 \leq k \leq \mathcal{L}} \sqrt{\sum_{\substack{1 < m \leq x , \\ p|m \Rightarrow x^{e^{-(k+1)}} < p \leq x^{e^{-k}}}} d_{\lceil q \rceil}(m) \Biggl( \frac{X}{m} \int_{m}^{m(1 + \frac{1}{X})} \E \Biggl| \sum_{\substack{x/t < n \leq x/m, \\ x^{e^{-(k+1)}} \text{-smooth}}} f(n) \Biggr|^{2q} dt \Biggr)^{1/q} } . \nonumber
\end{eqnarray}
The length of the sum over $n$ here is $\frac{x(t-m)}{mt} \leq \frac{x}{mX}$, so when $x/X \leq m \leq x$ there will be at most one term in the sum, and we simply have $\E \Biggl| \sum_{\substack{x/t < n \leq x/m, \\ x^{e^{-(k+1)}} \text{-smooth}}} f(n) \Biggr|^{2q} \leq 1$. When $x^{e^{-(k+1)}} < m < x/X$, we take a fairly crude approach and use the Cauchy--Schwarz inequality, obtaining that $\E \Biggl| \sum_{\substack{x/t < n \leq x/m, \\ x^{e^{-(k+1)}} \text{-smooth}}} f(n) \Biggr|^{2q}$ is at most
\begin{eqnarray}
\sqrt{\E \Biggl| \sum_{\substack{x/t < n \leq x/m, \\ x^{e^{-(k+1)}} \text{-smooth}}} f(n) \Biggr|^{2} \E \Biggl| \sum_{\substack{x/t < n \leq x/m, \\ x^{e^{-(k+1)}} \text{-smooth}}} f(n) \Biggr|^{2(2q-1)}} & \ll & \sqrt{\frac{x}{mX} \E \Biggl| \sum_{\substack{x/t < n \leq x/m, \\ x^{e^{-(k+1)}} \text{-smooth}}} f(n) \Biggr|^{2(2q-1)}} \nonumber \\
& \ll & \sqrt{\frac{x}{mX} (\frac{x}{m})^{2q-1} \log^{O(q^2)}x} = (\frac{x}{m})^{q} \frac{\log^{O(q^2)}x}{X^{1/2}} . \nonumber
\end{eqnarray}
Here the crude upper bound $(x/m)^{2q-1} \log^{O(q^2)}x$ for the $2(2q-1)$-th moment may be proved as in section \ref{seceasycases}.

Putting things together, we find that the second term in \eqref{uppersmoothingdisplay} is
$$ \ll \sum_{0 \leq k \leq \mathcal{L}} \sqrt{\frac{x \log^{O(q)}x}{X^{1/2q}} \sum_{\substack{1 < m < x/X , \\ p|m \Rightarrow x^{e^{-(k+1)}} < p \leq x^{e^{-k}}}} \frac{d_{\lceil q \rceil}(m)}{m} + \sum_{\substack{x/X \leq m \leq x , \\ p|m \Rightarrow x^{e^{-(k+1)}} < p \leq x^{e^{-k}}}} d_{\lceil q \rceil}(m) } . $$
To bound the first of these sums we use the simple estimate $\sum_{\substack{1 < m < x/X , \\ p|m \Rightarrow x^{e^{-(k+1)}} < p \leq x^{e^{-k}}}} \frac{d_{\lceil q \rceil}(m)}{m} \leq \prod_{x^{e^{-(k+1)}} < p \leq x^{e^{-k}}} (1 - \frac{1}{p})^{-\lceil q \rceil} = e^{O(q)}$. To bound the second sum, by sub-multiplicativity of $d_{\lceil q \rceil}(\cdot)$ we always have $d_{\lceil q \rceil}(m) \leq \lceil q \rceil^{\Omega(m)}$, and we note (to obtain good dependence on $k$) that if $m \geq x/X$ only has prime factors from the interval $(x^{e^{-(k+1)}}, x^{e^{-k}}]$, then we must have $\Omega(m) \geq e^{k}/2$, say. So we get that $\sum_{\substack{x/X \leq m \leq x , \\ p|m \Rightarrow x^{e^{-(k+1)}} < p \leq x^{e^{-k}}}} d_{\lceil q \rceil}(m)$ is at most
$$ 5^{-e^{k}/2} \sum_{\substack{x/X \leq m \leq x , \\ p|m \Rightarrow x^{e^{-(k+1)}} < p \leq x^{e^{-k}}}} (5\lceil q \rceil)^{\Omega(m)} \ll 5^{-e^{k}/2} \frac{e^{k} q x}{\log x} \prod_{  x^{e^{-(k+1)}} < p \leq x^{e^{-k}}} \left(1 - \frac{5 \lceil q \rceil}{p}\right)^{-1} \ll \frac{e^{O(q)} 2^{-e^{k}} x}{\log x} , $$
where the first inequality uses Number Theory Result 1. Recalling that we have $q \leq \log^{0.05}x$, and $\mathcal{L} = \lfloor (\log\log x)/10 \rfloor$, and $X = e^{\sqrt{\log x}}$, the second term in \eqref{uppersmoothingdisplay} is
$$ \leq \sum_{0 \leq k \leq \mathcal{L}} \sqrt{\frac{x \log^{O(q)}x}{X^{1/2q}} + e^{O(q)} 2^{-e^{k}} \frac{x}{\log x} } \leq e^{O(q)} \sqrt{\frac{x}{\log x}} , $$
which is an acceptable contribution for Proposition 1.

Turning to the remaining first sum in \eqref{uppersmoothingdisplay}, this is equal to
$$ \sum_{0 \leq k \leq \mathcal{L}} \vert\vert \int_{x^{e^{-(k+1)}}}^{x} \Biggl| \sum_{\substack{n \leq x/t, \\ x^{e^{-(k+1)}} \text{-smooth}}} f(n) \Biggr|^2 \sum_{\substack{t/(1+1/X) \leq m \leq t , \\ p|m \Rightarrow x^{e^{-(k+1)}} < p \leq x^{e^{-k}}}} \frac{X}{m} d_{\lceil q \rceil}(m) dt \vert\vert_{q}^{1/2} . $$
Now we set $u = u(k,t) := e^{k} (\log t)/\log x$, and notice that (by sub-multiplicativity) $d_{\lceil q \rceil}(m) \leq \lceil q \rceil^{\Omega(m)}$, and if $m \geq t/(1+1/X)$ only has prime factors from the interval $(x^{e^{-(k+1)}}, x^{e^{-k}}]$ then we must have $\Omega(m) \geq u - 1$. So using Number Theory Result 1 (whose conditions are satisfied since $X = e^{\sqrt{\log x}}$ isn't too large, and $k \leq (\log\log x)/10$) we get
\begin{eqnarray}
\sum_{\substack{t/(1+1/X) \leq m \leq t , \\ p|m \Rightarrow x^{e^{-(k+1)}} < p \leq x^{e^{-k}}}} \frac{X}{m} d_{\lceil q \rceil}(m) & \ll & \frac{X}{t} 5^{-u} \sum_{\substack{t/(1+1/X) \leq m \leq t , \\ p|m \Rightarrow x^{e^{-(k+1)}} < p \leq x^{e^{-k}}}} (5 \lceil q \rceil)^{\Omega(m)} \nonumber \\
& \ll & \frac{q e^k 5^{-u}}{\log x} \prod_{x^{e^{-(k+1)}} < p \leq x^{e^{-k}}} (1-\frac{5 \lceil q \rceil}{p})^{-1} \ll \frac{e^{O(q)}}{\log t} , \nonumber
\end{eqnarray}
provided $x$ is sufficiently large. Consequently, the first sum in \eqref{uppersmoothingdisplay} is
\begin{eqnarray}
& \leq & e^{O(q)} \sum_{0 \leq k \leq \mathcal{L}} \vert\vert \int_{x^{e^{-(k+1)}}}^{x} \Biggl| \sum_{\substack{n \leq x/t, \\ x^{e^{-(k+1)}} \text{-smooth}}} f(n) \Biggr|^2 \frac{dt}{\log t} \vert\vert_{q}^{1/2}  \nonumber \\
& = & e^{O(q)} \sum_{0 \leq k \leq \mathcal{L}} \sqrt{x} \vert\vert \int_{1}^{x^{1 - e^{-(k+1)}}} \Biggl| \sum_{\substack{n \leq z, \\ x^{e^{-(k+1)}} \text{-smooth}}} f(n) \Biggr|^2 \frac{dz}{z^2 \log(x/z)} \vert\vert_{q}^{1/2} , \nonumber
\end{eqnarray}
where the second line follows from making the substitution $z=x/t$.

To obtain a satisfactory dependence on $k$ in our final estimations, we now note that if $z \leq \sqrt{x}$ we have $\log(x/z) \gg \log x$, whereas if $\sqrt{x} < z \leq x^{1 - e^{-(k+1)}}$ we have $\log(x/z) \gg e^{-k} \log x$. Thus in any case we have $\log(x/z) \gg z^{-2k/\log x} \log x$. As discussed earlier, we also want to introduce a Rankin style shift, which we will achieve by adding a factor $(x/z)^{2q/\log x} = e^{O(q)} z^{-2q/\log x}$ into the integral. Inserting these estimates, we find the first sum in \eqref{uppersmoothingdisplay} is
$$ \leq \frac{\sqrt{x} e^{O(q)}}{\sqrt{\log x}} \sum_{0 \leq k \leq \mathcal{L}} \vert\vert \int_{1}^{x^{1 - e^{-(k+1)}}} \Biggl| \sum_{\substack{n \leq z, \\ x^{e^{-(k+1)}} \text{-smooth}}} f(n) \Biggr|^2 \frac{dz}{z^{2 + 2q/\log x - 2k/\log x}} \vert\vert_{q}^{1/2} . $$

Finally, using Harmonic Analysis Result 1 and then Minkowski's inequality, all of the above is
\begin{eqnarray}
& \leq & \sqrt{\frac{x}{\log x}} e^{O(q)} \sum_{0 \leq k \leq \mathcal{L}} \vert\vert \int_{-\infty}^{\infty} \frac{|F_k(1/2+q/\log x - k/\log x + it)|^2}{|1/2 + q/\log x -k/\log x + it|^2} dt \vert\vert_{q}^{1/2} \nonumber \\
& \leq & \sqrt{\frac{x}{\log x}} e^{O(q)} \sum_{0 \leq k \leq \mathcal{L}} \sqrt{\sum_{n \in \Z} \frac{1}{n^2 + 1} \vert\vert \int_{n-1/2}^{n+1/2} |F_k(1/2+\frac{q-k}{\log x} + it)|^2 dt \vert\vert_{q}} , \nonumber
\end{eqnarray}
where $F_k$ denotes the partial Euler product of $f(n)$ over $x^{e^{-(k+1)}}$-smooth numbers. In the Steinhaus case, since the law of the random function $f(n)$ is the same as the law of $f(n)n^{it}$ for any fixed $t \in \R$ we have
$$ \vert\vert \int_{n-1/2}^{n + 1/2} |F_k(1/2+ \frac{q-k}{\log x} + it)|^2 dt \vert\vert_{q} = \vert\vert \int_{-1/2}^{1/2} |F_k(1/2+\frac{q-k}{\log x} + it)|^2 dt \vert\vert_{q} \;\;\; \forall \; n . $$
Proposition 1 now follows on putting everything together.
\qed

\vspace{12pt}
The proof of Proposition 2, covering the Rademacher case, is very similar to the Steinhaus case. We use the Rademacher part of Probability Result 1, producing various terms $d_{2\lceil q \rceil - 1}(n)$ in place of $d_{\lceil q \rceil}(n)$, but this doesn't alter the analysis. The only non-trivial change comes at the very end of the proof, where (since it is no longer the case that the law of the random function $f(n)$ is the same as the law of $f(n)n^{it}$) we apply the bound
\begin{eqnarray}
&& \sqrt{\sum_{n \in \Z} \frac{1}{n^2 + 1} \vert\vert \int_{n-1/2}^{n+1/2} |F_k(1/2+\frac{q-k}{\log x} + it)|^2 dt \vert\vert_{q}} \nonumber \\
& \ll & \sqrt{\max_{N \in \Z} \frac{1}{(|N|+1)^{1/4}} \vert\vert \int_{N-1/2}^{N+1/2} |F_k(1/2+\frac{q-k}{\log x} + it)|^2 dt \vert\vert_{q} } . \nonumber
\end{eqnarray}
\qed

\subsection{Proof of Propositions 3 and 4}
We proceed somewhat similarly as in section 2.5 of Harper~\cite{harperrmflow}, or section 2.2 of Harper, Nikeghbali and Radziwi\l\l~\cite{hnr}.

Again we let $P(n)$ denote the largest prime factor of $n$, and we introduce an auxiliary Rademacher random variable $\epsilon$ that is independent of everything else. Then we find that
\begin{eqnarray}
\vert\vert \sum_{\substack{n \leq x, \\ P(n) > x^{3/4}}} f(n) \vert\vert_{2q} & = & \frac{1}{2} \vert\vert \sum_{\substack{n \leq x, \\ P(n) > x^{3/4}}} f(n) + \sum_{\substack{n \leq x, \\ P(n) \leq x^{3/4}}} f(n) + \sum_{\substack{n \leq x, \\ P(n) > x^{3/4}}} f(n) - \sum_{\substack{n \leq x, \\ P(n) \leq x^{3/4}}} f(n) \vert\vert_{2q} \nonumber \\
& \leq & \frac{1}{2} \Biggl( \vert\vert \sum_{\substack{n \leq x, \\ P(n) > x^{3/4}}} f(n) + \sum_{\substack{n \leq x, \\ P(n) \leq x^{3/4}}} f(n) \vert\vert_{2q} + \vert\vert \sum_{\substack{n \leq x, \\ P(n) > x^{3/4}}} f(n) - \sum_{\substack{n \leq x, \\ P(n) \leq x^{3/4}}} f(n) \vert\vert_{2q} \Biggr) \nonumber \\
& \leq & \vert\vert \epsilon \sum_{\substack{n \leq x, \\ P(n) > x^{3/4}}} f(n) + \sum_{\substack{n \leq x, \\ P(n) \leq x^{3/4}}} f(n) \vert\vert_{2q} = \vert\vert \sum_{n \leq x} f(n) \vert\vert_{2q} . \nonumber
\end{eqnarray}
Here the first inequality is Minkowski's inequality; the second is H\"{o}lder's inequality (with exponent $2q$) applied only to the averaging over $\epsilon$; and the final equality follows since the law of $\epsilon \sum_{\substack{n \leq x, \\ P(n) > x^{3/4}}} f(n) = \epsilon \sum_{x^{3/4} < p \leq x} f(p) \sum_{m \leq x/p} f(m)$ conditional on the values $(f(p))_{p \leq x^{3/4}}$ is the same as the law of $\sum_{\substack{n \leq x, \\ P(n) > x^{3/4}}} f(n)$.

Now in the decomposition $\sum_{\substack{n \leq x, \\ P(n) > x^{3/4}}} f(n) = \sum_{x^{3/4} < p \leq x} f(p) \sum_{m \leq x/p} f(m)$, the inner sums are determined by the values $(f(p))_{p \leq x^{3/4}}$ (and in fact by the values $(f(p))_{p \leq x^{1/4}}$), which are independent of the outer random variables $(f(p))_{x^{3/4} < p \leq x}$. So conditioning on the values $(f(p))_{p \leq x^{3/4}}$ determining the inner sums and applying Probability Result 2 with $a_p = \sum_{m \leq x/p} f(m)$, it follows that
$$ \vert\vert \sum_{\substack{n \leq x, \\ P(n) > x^{3/4}}} f(n) \vert\vert_{2q} \geq \vert\vert \sum_{x^{3/4} < p \leq x} \left|\sum_{m \leq x/p} f(m) \right|^2 \vert\vert_{q}^{1/2} \geq \frac{1}{\sqrt{\log x}} \vert\vert \sum_{x^{3/4} < p \leq x} \log p \left|\sum_{m \leq x/p} f(m) \right|^2 \vert\vert_{q}^{1/2} . $$

Next we want to replace the sum over $p$ by an integral average. We can rewrite
$$ \sum_{x^{3/4} < p \leq x} \log p \left|\sum_{m \leq x/p} f(m) \right|^2 = \sum_{r \leq x^{1/4}} \sum_{\frac{x}{r+1} < p \leq \frac{x}{r}} \log p \left|\sum_{m \leq r} f(m) \right|^2 , $$
and noting that $\frac{x}{r} - \frac{x}{r+1} = \frac{x}{r(r+1)} \gg (x/r)^{2/3}$ on our range of $r$, a Hoheisel-type Prime Number Theorem in short intervals (see e.g. Theorem 12.8 of Ivi\'c~\cite{ivic}) implies that
$$ \sum_{x^{3/4} < p \leq x} \log p \left|\sum_{m \leq x/p} f(m) \right|^2 \gg \sum_{r \leq x^{1/4}} \int_{x/(r+1)}^{x/r} 1 dt \left|\sum_{m \leq r} f(m) \right|^2 \geq \int_{x^{3/4}}^{x} \left|\sum_{m \leq x/t} f(m) \right|^2 dt . $$
Making a substitution $z = x/t$, we see this integral is the same as $x \int_{1}^{x^{1/4}} \left|\sum_{m \leq z} f(m) \right|^2 \frac{dz}{z^2}$. Checking back, this completes the proof of the first part of Proposition 3.

To deduce the second part of Proposition 3, we note that for any large $V$ and any $q \geq 1$ we have
\begin{eqnarray}
&& \vert\vert \int_{1}^{x^{1/4}} \left|\sum_{m \leq z} f(m) \right|^2 \frac{dz}{z^{2}} \vert\vert_{q} \geq \vert\vert \int_{1}^{x^{1/4}} \Biggl|\sum_{\substack{m \leq z, \\ x \text{-smooth}}} f(m) \Biggr|^2 \frac{dz}{z^{2 + 8Vq/\log x}} \vert\vert_{q} \nonumber \\
& \geq & \vert\vert \int_{1}^{\infty} \Biggl|\sum_{\substack{m \leq z, \\ x \text{-smooth}}} f(m) \Biggr|^2 \frac{dz}{z^{2 + 8Vq/\log x}} \vert\vert_{q} - \vert\vert \int_{x^{1/4}}^{\infty} \Biggl|\sum_{\substack{m \leq z, \\ x \text{-smooth}}} f(m) \Biggr|^2 \frac{dz}{z^{2 + 8Vq/\log x}} \vert\vert_{q} \nonumber \\
& \geq & \vert\vert \int_{1}^{\infty} \Biggl|\sum_{\substack{m \leq z, \\ x \text{-smooth}}} f(m) \Biggr|^2 \frac{dz}{z^{2 + 8Vq/\log x}} \vert\vert_{q} - \frac{1}{e^{Vq}} \vert\vert \int_{1}^{\infty} \Biggl|\sum_{\substack{m \leq z, \\ x \text{-smooth}}} f(m) \Biggr|^2 \frac{dz}{z^{2 + 4Vq/\log x}} \vert\vert_{q} . \nonumber
\end{eqnarray}
By Harmonic Analysis Result 1, provided that $V \leq (\log x)/q$ (so that $\frac{Vq}{\log x}$ is uniformly bounded) the first term here is $\gg \vert\vert \int_{-1/2}^{1/2} |F(1/2 + \frac{4Vq}{\log x} + it)|^2 dt \vert\vert_{q}$ and the subtracted second term is $\ll e^{-Vq} \vert\vert \int_{-\infty}^{\infty} \frac{|F(1/2 + \frac{2Vq}{\log x} + it)|^2}{|1/2 + \frac{2Vq}{\log x} + it|^2} dt \vert\vert_{q}$, which in the Steinhaus case is $\ll e^{-Vq} \vert\vert \int_{-1/2}^{1/2} |F(1/2 + \frac{2Vq}{\log x} + it)|^2 dt \vert\vert_{q} $ by ``translation invariance in law''. Putting everything together, this finishes the proof of Proposition 3.
\qed

\vspace{12pt}
The arguments in the Rademacher case are exactly the same until the final line, where we don't have ``translation invariance'' so we must upper bound $\vert\vert \int_{-\infty}^{\infty} \frac{|F(1/2 + \frac{2Vq}{\log x} + it)|^2}{|1/2 + \frac{2Vq}{\log x} + it|^2} dt \vert\vert_{q}$ by $\max_{N \in \Z} \frac{1}{(|N|+1)^{1/4}} \vert\vert \int_{N-1/2}^{N+1/2} |F(1/2 + \frac{2Vq}{\log x} + it)|^2 dt \vert\vert_{q}$, say.
\qed

\section{Proofs of the upper bounds in Theorems 1 and 2}\label{secmainupper}
In view of Proposition \ref{propstupper}, the key to obtaining the upper bound in Theorem 1 will lie in proving the following. Recall here that $F_{k}(s)$ denotes the partial Euler product of $f(n)$ over $x^{e^{-(k+1)}}$-smooth numbers, and in the special case where $k=-1$ we usually write $F(s)$ (rather than $F_{-1}(s)$) for the partial Euler product over $x$-smooth numbers.

\begin{keyprop1}
Let $f(n)$ be a Steinhaus random multiplicative function. For all large $x$, and uniformly for $1 \leq q \leq \log^{100}x$ (say) and $-1 \leq k \leq \mathcal{L} = \lfloor (\log\log x)/10 \rfloor$ and $- \frac{e^k}{\log x} \leq \sigma \leq \frac{1}{100\log(2q)}$ (say), we have
$$ \E( \int_{-1/2}^{1/2} |F_{k}(1/2 + \sigma + it)|^2 dt )^q \ll \frac{e^{O(q^2)}}{\log^{q-1}x} \Biggl( \frac{\log x}{\log 2q} \Biggr)^{q^2} \min\Biggl\{\frac{1}{e^{k+1}}, \frac{1}{|\sigma| \log x} \Biggr\}^{q^2 - q + 1} . $$
\end{keyprop1}

Key Proposition 1 is actually much more general, in terms of the allowed range of $q$ and $\sigma$, than we immediately need (and the proof would let us extend the range of $q$ quite a lot further if we wished, all it really requires is something like $\frac{e^k}{\log x} \leq \frac{1}{100\log(2q)}$). The increased generality will be useful in section \ref{secmainlower}, where Key Proposition 1 will play an auxiliary role, and also in clarifying the essential features of the proof.

\begin{proof}[Proof of the upper bound in Theorem 1, assuming Key Proposition 1]
In view of the discussion in section \ref{seceasycases}, it will suffice to prove the Theorem 1 upper bound for $1 \leq q \leq \log\log x$. And to do that, in view of Proposition \ref{propstupper} it will suffice to show that
$$ \sum_{0 \leq k \leq \mathcal{L}} \vert\vert \int_{-1/2}^{1/2} |F_k(1/2 + \frac{q-k}{\log x} + it)|^2 dt \vert\vert_{q}^{1/2} \leq e^{-(q/2)\log q - (q/2)\log\log(2q) + O(q)} \log^{q/2 - 1/2 + 1/2q}x . $$

Applying Key Proposition 1 with $\sigma = \frac{q-k}{\log x}$ (which is indeed $\leq \frac{1}{100\log(2q)}$ on our range of $q$), we find the left hand side is
\begin{eqnarray}
& \leq & \sum_{0 \leq k \leq \mathcal{L}} \left( \frac{e^{O(q^2)}}{\log^{q-1}x} \Biggl( \frac{\log x}{\log 2q} \Biggr)^{q^2} \min\Biggl\{\frac{1}{e^{k+1}}, \frac{1}{|q-k|} \Biggr\}^{q^2 - q + 1} \right)^{1/2q} \nonumber \\
& = & \sum_{0 \leq k \leq \mathcal{L}} \left( e^{O(q^2)} \Biggl( \frac{\log x}{\log 2q} \min\Biggl\{\frac{1}{e^{k+1}}, \frac{1}{q} \Biggr\} \Biggr)^{q^2 - q + 1} \right)^{1/2q} . \nonumber
\end{eqnarray}
It is easy to see that this satisfies our desired bound.
\end{proof}

For Theorem 2, we need a Rademacher analogue of the above.

\begin{keyprop2}
Let $f(n)$ be a Rademacher random multiplicative function. For all large $x$, and uniformly for $1 \leq q \leq \log^{100}x$ (say) and $-1 \leq k \leq \mathcal{L} = \lfloor (\log\log x)/10 \rfloor$ and $- \frac{e^k}{\log x} \leq \sigma \leq \frac{1}{100\log(2q)}$ (say), we have
\begin{eqnarray}
\E( \int_{-1/2}^{1/2} |F_{k}(\frac{1}{2} + \sigma + it)|^2 dt )^q & \ll & \frac{e^{O(q^2)}}{\log^{q}x} (1 + \min\{\log\log x, \frac{1}{|q-q_0|}\}) \Biggl( \frac{\log x}{\log 2q} \Biggr)^{\max\{2q^2 - q, q^2 +1\}} \cdot \nonumber \\
&& \cdot \min\Biggl\{\frac{1}{e^{k+1}}, \frac{1}{|\sigma| \log x} \Biggr\}^{\max\{2q^2 - 2q, q^2 - q + 1\}} , \nonumber
\end{eqnarray}
where $q_0 = (1 + \sqrt{5})/2$.

Furthermore, for any $|N| \geq 1$ we have
\begin{eqnarray}
\E( \int_{N-\frac{1}{2}}^{N+\frac{1}{2}} |F_{k}(\frac{1}{2} + \sigma + it)|^2 dt )^q & \ll & \min\Biggl\{|N|^{\frac{1}{100}}, \frac{\log x}{e^{k+1} \log 2q}, \frac{1}{|\sigma| \log 2q}\Biggr\}^{q(q+1)} \cdot \nonumber \\
&& \cdot \frac{e^{O(q^2)}}{\log^{q-1}x} \Biggl( \frac{\log x}{\log 2q} \Biggr)^{q^2} \min\Biggl\{\frac{1}{e^{k+1}}, \frac{1}{|\sigma| \log x} \Biggr\}^{q^2 - q + 1} . \nonumber
\end{eqnarray}
\end{keyprop2}

\begin{proof}[Proof of the upper bound in Theorem 2, assuming Key Proposition 2]
Similarly as in the Steinhaus case, in view of Proposition \ref{propradupper} and the discussion in section \ref{seceasycases} it will suffice to show that for all $1 \leq q \leq \log\log x$, we have
\begin{eqnarray}
&& \sum_{0 \leq k \leq \mathcal{L}} \max_{N \in \Z} \frac{1}{(|N|+1)^{1/8}} \vert\vert \int_{N-1/2}^{N+1/2} |F_k(1/2 + \frac{q-k}{\log x} + it)|^2 dt \vert\vert_{q}^{1/2} \nonumber \\
& \leq & e^{-q\log q - q\log\log(2q) + O(q)} (1 + \min\{\log\log x, \frac{1}{|q-q_0|}\})^{1/2q} (\log x)^{\max\{q - 1, q/2 - 1/2 + 1/2q\}} . \nonumber
\end{eqnarray}

We apply Key Proposition 2 with $\sigma = \frac{q-k}{\log x}$ (which is indeed $\leq \frac{1}{100\log(2q)}$ on our range of $q$). When $q \leq 15$, say, we have $\min\{|N|^{\frac{1}{100}}, \frac{\log x}{e^{k+1} \log 2q}, \frac{1}{|\sigma| \log 2q}\}^{q(q+1)} \leq |N|^{q/5} \leq |N|^{2q/8}$, so (on taking $2q$-th roots in Key Proposition 2 and then multiplying by the prefactor $1/(|N|+1)^{1/8}$) we see the contribution from $|N| \geq 1$ to $\max_{N \in \Z}$ will never exceed the contribution from the $N=0$ term. So overall, when $1 \leq q \leq 15$ the left hand side will be
\begin{eqnarray}
& \leq & \sum_{0 \leq k \leq \mathcal{L}} \left( \frac{e^{O(q^2)}}{\log^{q}x} \min\{\log\log x, \frac{1}{|q-q_0|}\} \Biggl( \frac{\log x}{\log 2q} \Biggr)^{\max\{2q^2-q,q^2+1\}} \Biggl( \frac{1}{e^{k+1}} \Biggr)^{\max\{2q^2 - 2q, q^2 - q + 1\}} \right)^{\frac{1}{2q}} \nonumber \\
& \ll & \sum_{0 \leq k \leq \mathcal{L}} \left( \min\{\log\log x, \frac{1}{|q-q_0|}\} \Biggl( \frac{\log x}{e^{k+1}} \Biggr)^{\max\{2q(q-1),q^2 - q + 1\}} \right)^{1/2q} . \nonumber
\end{eqnarray}
This certainly gives our desired bound for $1 \leq q \leq 15$.

When $15 \leq q \leq \log\log x$, we note first that $\max\{2q^2 - 2q, q^2 - q + 1\} = 2q^2 - 2q$. (In fact this is true as soon as $q \geq q_0$.) So using the bound $\min\{|N|^{\frac{1}{100}}, \frac{\log x}{e^{k+1} \log 2q}, \frac{1}{|\sigma| \log 2q}\}^{q(q+1)} \leq |N|^{(2q+1)/100} \min\{\frac{\log x}{e^{k+1} \log 2q}, \frac{1}{|\sigma| \log 2q}\}^{q^2 - q - 1}$, we again find that the contribution from $|N| \geq 1$ to $\max_{N \in \Z}$ will never exceed the contribution from the $N=0$ term. Overall, in this case we get a bound
$$ \ll \sum_{0 \leq k \leq \mathcal{L}} \left( e^{O(q^2)} \Biggl( \frac{\log x}{\log 2q} \min\Biggl\{\frac{1}{e^{k+1}}, \frac{1}{q} \Biggr\} \Biggr)^{2q^2 - 2q} \right)^{1/2q} \ll \left( e^{O(q^2)} \Biggl( \frac{\log x}{q \log 2q} \Biggr)^{2q^2 - 2q} \right)^{1/2q} , $$
as desired.
\end{proof}

We shall prove Key Propositions 1 and 2 in several steps over the course of this section. For convenience in the writing we set $\mathcal{X} := \min\{\frac{\log x}{e^{k+1}}, \frac{1}{|\sigma|}\}$, and note that under our hypotheses this is always $\geq 100\log(2q)$. The point of this definition is that the contribution from primes $p > e^{\mathcal{X}}$ in our Euler products will ultimately contribute only to the $e^{O(q^2)}$ term. 

\subsection{Preliminary manoeuvres}
We begin with a few manipulations to discretise and set up the problem, in both the Steinhaus and Rademacher cases. For any $q \geq 1$, we have
\begin{eqnarray}
\vert\vert \int_{-1/2}^{1/2} |F_k(1/2+\sigma + it)|^2 dt \vert\vert_{q} & \leq & \vert\vert \int_{-\frac{1}{2\mathcal{X}}}^{\frac{1}{2\mathcal{X}}} \sum_{|n| \leq \mathcal{X}/2 + 1} |F_k(1/2+\sigma + i(\frac{n}{\mathcal{X}} + t))|^2 dt \vert\vert_{q} \nonumber \\
& = & \frac{1}{\mathcal{X}} \vert\vert \int_{-\frac{1}{2\mathcal{X}}}^{\frac{1}{2\mathcal{X}}} \mathcal{X} \sum_{|n| \leq \mathcal{X}/2 + 1} |F_k(1/2+\sigma + i(\frac{n}{\mathcal{X}} + t))|^2 dt \vert\vert_{q} . \nonumber
\end{eqnarray}
Applying H\"{o}lder's inequality with exponent $q$ to the normalised integral $\int_{-\frac{1}{2\mathcal{X}}}^{\frac{1}{2\mathcal{X}}} \mathcal{X} dt$, we see the right hand side is
$$ \leq \frac{1}{\mathcal{X}} \left( \int_{-\frac{1}{2\mathcal{X}}}^{\frac{1}{2\mathcal{X}}} \mathcal{X} \; \E \left(\sum_{|n| \leq \mathcal{X}/2 + 1} |F_k(1/2+\sigma + i(\frac{n}{\mathcal{X}} + t))|^2 \right)^q dt \right)^{1/q} . $$

In the Steinhaus case, where $|F_k(1/2+\sigma + i(\frac{n}{\mathcal{X}} + t))|^2$ has the same distribution for any given shift $t$, we can simplify the above to give the bound
$$ \vert\vert \int_{-1/2}^{1/2} |F_k(1/2+\sigma + it)|^2 dt \vert\vert_{q} \leq \frac{1}{\mathcal{X}} \vert\vert \sum_{|n| \leq \mathcal{X}/2 + 1} |F_k(1/2+\sigma + \frac{in}{\mathcal{X}} )|^2 \vert\vert_{q} , $$
so we have indeed passed to studying a discrete sum rather than an integral. Finally, we rewrite the right hand side as
\begin{eqnarray}\label{upperq-1}
&& \frac{1}{\mathcal{X}} \left(\E \sum_{|n| \leq \frac{\mathcal{X}}{2} + 1} |F_k(\frac{1}{2}+\sigma + \frac{in}{\mathcal{X}})|^2 \left( \sum_{|m| \leq \frac{\mathcal{X}}{2} + 1} |F_k(\frac{1}{2}+\sigma + \frac{im}{\mathcal{X}})|^2 \right)^{q-1} \right)^{1/q} \nonumber \\
& \ll & \frac{1}{\mathcal{X}} \left( \mathcal{X} \; \E |F_k(\frac{1}{2}+\sigma)|^2 \left( \sum_{|m| \leq \mathcal{X}} |F_k(\frac{1}{2}+\sigma + \frac{im}{\mathcal{X}})|^2 \right)^{q-1} \right)^{1/q} ,
\end{eqnarray}
where the inequality again uses the distributional ``translation invariance'' (shifting $n$ to zero in the outer sum, and replacing $m$ by $m-n$ in the second sum).

In the case of Rademacher $f(n)$, if we mimic the above calculations we obtain that $\vert\vert \int_{N-1/2}^{N+1/2} |F_k(\frac{1}{2}+\sigma + it)|^2 dt \vert\vert_{q}$ is
\begin{eqnarray}\label{radupperq-1}
& \leq & \frac{1}{\mathcal{X}} \left(\int_{-\frac{1}{2\mathcal{X}}}^{\frac{1}{2\mathcal{X}}} \mathcal{X} \; \E \sum_{|n| \leq \frac{\mathcal{X}}{2} + 1} |F_k(\frac{1}{2}+\sigma + i(\frac{n}{\mathcal{X}} + N + t))|^2 \cdot \right. \nonumber \\
&& \left. \cdot \left( \sum_{|m| \leq \frac{\mathcal{X}}{2} + 1} |F_k(\frac{1}{2}+\sigma + i(\frac{m}{\mathcal{X}} + N + t))|^2 \right)^{q-1} dt \right)^{1/q} .
\end{eqnarray}

\subsection{Proof of Key Proposition 1, for $q \geq 2$}
When $q \geq 2$, we are helped by the fact that we can use H\"{o}lder's inequality again to analyse $\left( \sum_{|m| \leq \mathcal{X}} |F_k(\frac{1}{2}+\sigma + \frac{im}{\mathcal{X}})|^2 \right)^{q-1}$ in \eqref{upperq-1}. If we let $\mu := \sum_{|m| \leq \mathcal{X}} \frac{1}{(|m|+1)^2}$, so that $\mu \asymp 1$, then first we have
$$ \left( \sum_{|m| \leq \mathcal{X}} |F_k(\frac{1}{2}+\sigma + \frac{im}{\mathcal{X}})|^2 \right)^{q-1} = \mu^{q-1} \left( \frac{1}{\mu} \sum_{|m| \leq \mathcal{X}} \frac{1}{(|m|+1)^2} (|m|+1)^2 |F_k(\frac{1}{2}+\sigma + \frac{im}{\mathcal{X}})|^2 \right)^{q-1} . $$
Using H\"{o}lder's inequality with exponent $q-1$, we deduce $\left( \sum_{|m| \leq \mathcal{X}} |F_k(\frac{1}{2}+\sigma + \frac{im}{\mathcal{X}})|^2 \right)^{q-1}$ is
\begin{eqnarray}\label{holderbigq-1}
& \leq & \mu^{q-1} \cdot \frac{1}{\mu} \sum_{|m| \leq \mathcal{X}} \frac{1}{(|m|+1)^2} \left((|m|+1)^2 |F_k(1/2+\sigma + \frac{im}{\mathcal{X}})|^2 \right)^{q-1} \nonumber \\
& = & e^{O(q)} \sum_{|m| \leq \mathcal{X}} \frac{1}{(|m|+1)^2} (|m|+1)^{2(q-1)} |F_k(1/2+\sigma + \frac{im}{\mathcal{X}})|^{2(q-1)} .
\end{eqnarray}
We remark that the choice of weights $1/(|m|+1)^2$ that we introduced is fairly arbitrary. The key point is that we expect, in \eqref{upperq-1}, that the only significant contribution should come from small $m$ (for which $|F_k(\frac{1}{2}+\sigma + \frac{im}{\mathcal{X}})|^2$ will be highly correlated with the outer term $|F_k(\frac{1}{2}+\sigma)|^2$), so we don't want to pick up a factor like $\mathcal{X}$ (inefficiently reflecting the total length of the sum) in our application of H\"{o}lder's inequality.

In view of the above computation, to bound the right hand side of \eqref{upperq-1} when $q \geq 2$ we need to bound terms of the form
$$ (|m|+1)^{2(q-1)} \E |F_k(1/2+\sigma)|^{2} |F_k(1/2+\sigma + \frac{im}{\mathcal{X}})|^{2(q-1)} . $$
Recall that $- \frac{e^k}{\log x} \leq \sigma \leq \frac{1}{100\log(2q)}$ here. Inserting the definition of $F_k(s)$, and using a trivial bound $e^{O(\sum_{p \leq 100q^2} q/p^{1/2+\sigma})} = e^{O(\sum_{p \leq 100q^2} q/\sqrt{p})} = e^{O(q^{2}/\log q)}$ for the parts of the Euler products over primes $\leq 100q^2$, this is
$$ e^{O(q^{2}/\log q)} (|m|+1)^{2(q-1)} \E \prod_{100q^2 < p \leq x^{e^{-(k+1)}}} \left|1 - \frac{f(p)}{p^{1/2+\sigma}}\right|^{-2} \left|1 - \frac{f(p)}{p^{1/2+\sigma + im/\mathcal{X}}}\right|^{-2(q-1)} . $$
Now if $\frac{e^{k+1}}{\log x} \leq \sigma \leq \frac{1}{100\log(2q)}$, (and so $\mathcal{X} := \min\{\frac{\log x}{e^{k+1}}, \frac{1}{|\sigma|}\} = \frac{1}{\sigma}$), then the first part of Euler Product Result 1 implies that the expectation of the part of the Euler product over primes $e^{1/\sigma} < p \leq x^{e^{-(k+1)}}$ is equal to $\exp\{O(\sum_{e^{1/\sigma} < p \leq x^{e^{-(k+1)}}} \frac{q^2}{p^{1+2\sigma}} + \frac{q^3}{e^{1/2\sigma}})\}$, which is all $e^{O(q^2)}$. Using this fact, as well as the independence of $f(p)$ for different primes $p$, we find the above is always equal to
$$ e^{O(q^{2})} (|m|+1)^{2(q-1)} \E \prod_{100q^2 < p \leq e^{\mathcal{X}}} \left|1 - \frac{f(p)}{p^{1/2+\sigma}}\right|^{-2} \left|1 - \frac{f(p)}{p^{1/2+\sigma + im/\mathcal{X}}}\right|^{-2(q-1)} . $$
Notice that our size assumptions on $\sigma, k, q$ guarantee that $e^{\mathcal{X}}$ is larger than $100q^2$. Finally, the second part of Euler Product Result 1 implies this is all equal to
$$ e^{O(q^{2})} (|m|+1)^{2(q-1)} \left(\frac{\mathcal{X}}{\log q} \right)^{1 + (q-1)^2} \left(1 + \frac{\mathcal{X}}{(|m|+1)\log q} \right)^{2(q-1)} = e^{O(q^{2})} \left(\frac{\mathcal{X}}{\log q} \right)^{q^2} . $$

Putting this together with \eqref{holderbigq-1} and \eqref{upperq-1}, we find $\vert\vert \int_{-1/2}^{1/2} |F_k(1/2+\sigma + it)|^2 dt \vert\vert_{q}$ is
$$ \ll \frac{1}{\mathcal{X}} \left( \mathcal{X} \cdot e^{O(q^2)} \sum_{|m| \leq \mathcal{X}} \frac{1}{(|m|+1)^2} \left(\frac{\mathcal{X}}{\log q} \right)^{q^2} \right)^{1/q} = \frac{1}{\mathcal{X}} \left( \mathcal{X} \cdot e^{O(q^2)} \left(\frac{\mathcal{X}}{\log q} \right)^{q^2} \right)^{1/q} . $$
Raising everything to the power $q$ and inserting the fact that $\mathcal{X} = \min\{\frac{\log x}{e^{k+1}}, \frac{1}{|\sigma|}\}$, this gives the statement of Key Proposition 1.
\qed

\subsection{Proof of Key Proposition 2, for $q \geq 2$}
We begin with the second part of Key Proposition 2, where $|N| \geq 1$. Then similarly as in the deduction of \eqref{holderbigq-1} in the Steinhaus case, for any $|n| \leq \frac{\mathcal{X}}{2} + 1$ and any $|t| \leq 1/(2\mathcal{X})$ we can use H\"{o}lder's inequality to show
\begin{eqnarray}
&& \left( \sum_{|m| \leq \frac{\mathcal{X}}{2} + 1} |F_k(1/2+\sigma + i(\frac{m}{\mathcal{X}} + N + t))|^2 \right)^{q-1} \nonumber \\
& \leq & e^{O(q)} \sum_{|m| \leq \frac{\mathcal{X}}{2} + 1 } \frac{1}{(|m-n|+1)^2} (|m-n|+1)^{2(q-1)} |F_k(1/2+\sigma + i(\frac{m}{\mathcal{X}} + N + t))|^{2(q-1)} . \nonumber
\end{eqnarray}

So to bound the right hand side of \eqref{radupperq-1}, we need to bound terms of the form 
$$ (|m-n|+1)^{2(q-1)} \E |F_k(\frac{1}{2}+\sigma + i(\frac{n}{\mathcal{X}} + N + t))|^{2} |F_k(\frac{1}{2}+\sigma + i(\frac{m}{\mathcal{X}} + N + t))|^{2(q-1)} . $$
As in the Steinhaus case, the contribution from primes $p \leq 100q^2$ to this expectation is trivially $e^{O(q^{2}/\log q)}$. Using the first part of Euler Product Result 2, the contribution from primes $e^{\mathcal{X}} < p \leq x^{e^{-(k+1)}}$ is $e^{O(q^2)}$, and overall (noting that the imaginary shifts $\frac{n}{\mathcal{X}} + N + t, \frac{m}{\mathcal{X}} + N + t$ are $\gg |N| \gg 1$ and also $\ll |N|$) the above expression is at most
$$ e^{O(q^2)} \min\{\frac{\mathcal{X}}{\log q}, |N|^{\frac{1}{100}}\}^{(q-1)^2 + 3(q-1)} (|m-n|+1)^{2(q-1)} \left(\frac{\mathcal{X}}{\log q} \right)^{1 + (q-1)^2} \left(\frac{\mathcal{X}}{(|m-n|+1)\log q} \right)^{2(q-1)} . $$
Apart from the factor $\min\{\frac{\mathcal{X}}{\log q}, |N|^{\frac{1}{100}}\}^{(q-1)(q+2)}$, this is precisely analogous to the estimate we had in the Steinhaus case, so when $|N| \geq 1$ we indeed get the same bound as in the Steinhaus case apart from a multiplier $\min\{\frac{\mathcal{X}}{\log q}, |N|^{\frac{1}{100}}\}^{(q-1)(q+2)} \leq \min\{\frac{\mathcal{X}}{\log q}, |N|^{\frac{1}{100}}\}^{q(q+1)} = \min\{\frac{\log x}{e^{k+1} \log q}, \frac{1}{|\sigma| \log q}, |N|^{1/100}\}^{q(q+1)}$.

It remains to address the first part of Key Proposition 2, where $N=0$. In this case, when $q \geq 2$ we expect the main contribution to \eqref{radupperq-1} to come from terms with $n,m \approx 0$, so rather than splitting up the sum over $m$ according to the size of $|m-n|$ we shall just split it up according to the size of $m$. Proceeding in this way, using H\"{o}lder's inequality as in the Steinhaus case we find for any $|t| \leq 1/(2\mathcal{X})$ that
\begin{eqnarray}\label{radholderbigq-1}
&& \sum_{|n| \leq \frac{\mathcal{X}}{2} + 1} |F_k(\frac{1}{2}+\sigma + i(\frac{n}{\mathcal{X}} + t))|^2 \left( \sum_{|m| \leq \frac{\mathcal{X}}{2} + 1} |F_k(1/2+\sigma + i(\frac{m}{\mathcal{X}} + t))|^2 \right)^{q-1} \nonumber \\
& \leq & e^{O(q)} \sum_{|n| \leq \frac{\mathcal{X}}{2} + 1} |F_k(\frac{1}{2}+\sigma + i(\frac{n}{\mathcal{X}} + t))|^2 \cdot \nonumber \\
&& \cdot \sum_{|m| \leq \frac{\mathcal{X}}{2} + 1} \frac{1}{(|m|+1)^2} (|m|+1)^{2(q-1)} |F_k(1/2+\sigma + i(\frac{m}{\mathcal{X}} + t))|^{2(q-1)} .
\end{eqnarray}

So to bound the right hand side of \eqref{radupperq-1}, we again need to bound terms of the form $ (|m|+1)^{2(q-1)} \E |F_k(1/2+\sigma + i(\frac{n}{\mathcal{X}} + t))|^{2} |F_k(1/2+\sigma + i(\frac{m}{\mathcal{X}} + t))|^{2(q-1)}$. Using the second part of Euler Product Result 2, this is
\begin{eqnarray}
&& e^{O(q^2)} (|m|+1)^{2(q-1)} \left(1 + \frac{\mathcal{X}}{(|m|+1)\log q} \right)^{(q-1)^2 - (q-1)} \left(\frac{\mathcal{X}}{\log q} \right)^{1 + (q-1)^2} \nonumber \\
&& \cdot \left( \left(1 + \frac{\mathcal{X}}{(|m-n|+1)\log q} \right) \left(1 + \frac{\mathcal{X}}{(|m+n|+1)\log q} \right) \right)^{2(q-1)} . \nonumber
\end{eqnarray}
Now depending on the signs of $m,n$, one of the terms $|m-n|, |m+n|$ will be equal to $||m|-|n||$ and the other will equal $|m| + |n| \geq |m|$. So the above is always
$$ \leq e^{O(q^2)} (|m|+1)^{2(q-1)} \left(1 + \frac{\mathcal{X}}{(|m|+1)\log q} \right)^{q(q-1)} \left(\frac{\mathcal{X}}{\log q} \right)^{1 + (q-1)^2} \left(\frac{\mathcal{X}}{(||m|-|n||+1)\log q} \right)^{2(q-1)} . $$
Putting this together with \eqref{radholderbigq-1} and \eqref{radupperq-1}, if we first perform the sum over $|n| \leq \mathcal{X}/2 + 1$ we get that $\vert\vert \int_{-1/2}^{1/2} |F_k(1/2+\sigma + it)|^2 dt \vert\vert_{q}$ is at most
$$ \frac{e^{O(q)}}{\mathcal{X}} \left( \sum_{|m| \leq \frac{\mathcal{X}}{2} + 1} \frac{(|m|+1)^{2(q-1)}}{(|m|+1)^2} \left(1 + \frac{\mathcal{X}}{(|m|+1)\log q} \right)^{q(q-1)} \left(\frac{\mathcal{X}}{\log q} \right)^{1 + (q-1)^2 + 2(q-1)} \right)^{1/q} . $$
Finally performing the sum over $m$, the dominant contribution comes from small terms (note that terms with $|m| > \frac{\mathcal{X}}{\log q}$ contribute at most $e^{O(q^2)} \left(\frac{\mathcal{X}}{\log q}\right)^{2(q-1)} \left(\frac{\mathcal{X}}{\log q} \right)^{1 + (q-1)^2 + 2(q-1)}$ inside the bracket), and gives us a bound $\ll \frac{e^{O(q)}}{\mathcal{X}} \left(\frac{\mathcal{X}}{\log q}\right)^{(2q^2 - q)/q}$. Raising everything to the power $q \geq 2$ and inserting the fact that $\mathcal{X} = \min\{\frac{\log x}{e^{k+1}}, \frac{1}{|\sigma|}\}$, this gives the bound claimed in Key Proposition 2.
\qed

\subsection{Proof of Key Proposition 1, for $1 < q < 2$}
When $1 < q < 2$, it is not immediately obvious how to analyse the term $\E |F_k(\frac{1}{2}+\sigma)|^2 \left( \sum_{|m| \leq \mathcal{X}} |F_k(\frac{1}{2}+\sigma + \frac{im}{\mathcal{X}})|^2 \right)^{q-1}$ in \eqref{upperq-1}. We may begin by letting $C = C(q) = e^{1/(q-1)}$, and noting that $ \E |F_k(\frac{1}{2}+\sigma)|^2 \left( \sum_{|m| \leq \mathcal{X}} |F_k(\frac{1}{2}+\sigma + \frac{im}{\mathcal{X}})|^2 \right)^{q-1}$ is
$$ \leq \E |F_k(\frac{1}{2}+\sigma)|^2 \left( \sum_{d \leq (q-1)\log\mathcal{X} + 1} \sum_{C^{d-1} \leq |m| \leq C^{d}} |F_k(\frac{1}{2}+\sigma + \frac{im}{\mathcal{X}})|^2 \right)^{q-1} . $$
Here we adopt the convention that the term $m=0$ is included in $\sum_{C^{d-1} \leq |m| \leq C^{d}}$ when $d=1$, and that any terms with $|m| > \mathcal{X}$ are omitted from all sums (so the imaginary shift in the second copy of $F_k$ always has size $|m/\mathcal{X}| \leq 1$). The motivation for splitting things up like this is that we expect our estimates for all terms with $\sum_{C^{d-1} \leq |m| \leq C^{d}}$ to be roughly the same, up to a factor $C^{O(1)}$. And, when everything is raised to the power $q-1$, this factor simply becomes a constant multiplier. Next, if we let $D = D(q) \in \N$ be a parameter, to be fixed later, we can split things up further and find
\begin{eqnarray}\label{maxdsumm}
&& \E |F_k(\frac{1}{2}+\sigma)|^2 \left( \sum_{|m| \leq \mathcal{X}} |F_k(\frac{1}{2}+\sigma + \frac{im}{\mathcal{X}})|^2 \right)^{q-1} \nonumber \\
& \leq & \sum_{r \leq \frac{(q-1)\log\mathcal{X}}{D} + 1} \E |F_k(\frac{1}{2}+\sigma)|^2 \left( \sum_{(r-1)D < d \leq rD} \sum_{C^{d-1} \leq |m| \leq C^{d}} |F_k(\frac{1}{2}+\sigma + \frac{im}{\mathcal{X}})|^2 \right)^{q-1} \nonumber \\
& \ll & D^{q-1} \sum_{r \leq \frac{(q-1)\log\mathcal{X}}{D} + 1} \E \max_{(r-1)D < d \leq rD} |F_k(\frac{1}{2}+\sigma)|^2 \left( \sum_{C^{d-1} \leq |m| \leq C^{d}} |F_k(\frac{1}{2}+\sigma + \frac{im}{\mathcal{X}})|^2 \right)^{q-1} .
\end{eqnarray}
Notice that we may further assume that all terms $d$ for which $C^{d-1} > \mathcal{X}$ are omitted here, since for those the sum over $m$ is empty (by our earlier convention).

Now in the sum over $m$, we expect (thinking about Euler Product Result 1) that the part of the Euler product $F_k(\frac{1}{2}+\sigma + \frac{im}{\mathcal{X}})$ on primes $\leq e^{\mathcal{X}/C^d}$ will be roughly the same size for all $C^{d-1} \leq |m| \leq C^{d}$, and indeed roughly the same size as the corresponding part of $F_k(\frac{1}{2}+\sigma)$. To simplify our writing about this, for each $d \geq 1$ and $|m| \leq \mathcal{X}$ let us set 
$$ G_{d}(m) := \prod_{p \leq e^{\mathcal{X}/C^d}}  \left|1 - \frac{f(p)}{p^{1/2+\sigma+im/\mathcal{X}}}\right|^{-2} , \;\;\;\;\; \text{and} \;\;\;\;\; H_{d}(m) := \frac{|F_k(\frac{1}{2}+\sigma + \frac{im}{\mathcal{X}})|^2}{G_{d}(m)} . $$
(These quantities of course depend on $x,k,\sigma$ as well, but we suppress that in our notation.) We will also set $G_{d} := G_{d}(0)$ and $H_{d} := H_{d}(0)$. Then the expectation in \eqref{maxdsumm} may be written as
$$ \E \max_{(r-1)D < d \leq rD} G_{d} H_{d} \left( \sum_{C^{d-1} \leq |m| \leq C^{d}} G_{d}(m) H_{d}(m) \right)^{q-1} . $$
We want to apply H\"{o}lder's inequality to this expectation, in such a way that the bracketed sum is raised to the power $1/(q-1)$, and so we can connect up the expectation with the terms inside. Prior to doing this, we rewrite the expectation again as
$$ \E \max_{(r-1)D < d \leq rD} G_{d}^{1-(q-1)^2} H_{d}^{1-(q-1)} C^{2(q-1)(2-q)d} \left( \frac{G_{d}^{q-1} H_{d}}{C^{2(2-q)d}} \sum_{C^{d-1} \leq |m| \leq C^{d}} G_{d}(m) H_{d}(m) \right)^{q-1} . $$
Simplifying the various exponents, this is all
$$ \leq \E \left( \max_{(r-1)D < d \leq rD} G_{d}^{q(2-q)} H_{d}^{2-q} C^{2(q-1)(2-q)d} \right) \cdot \left( \max_{(r-1)D < d \leq rD} \frac{G_{d}^{q-1} H_{d}}{C^{2(2-q)d}} \sum_{C^{d-1} \leq |m| \leq C^{d}} G_{d}(m) H_{d}(m) \right)^{q-1} , $$
and now using H\"{o}lder's inequality with exponents $1/(2-q)$ and $1/(q-1)$, we get a bound
\begin{eqnarray}
& \leq & \left( \E \max_{(r-1)D < d \leq rD} G_{d}^{q} H_{d} C^{2(q-1)d} \right)^{2-q} \cdot \left( \E \max_{(r-1)D < d \leq rD} \frac{G_{d}^{q-1} H_{d}}{C^{2(2-q)d}} \sum_{C^{d-1} \leq |m| \leq C^{d}} G_{d}(m) H_{d}(m) \right)^{q-1} \nonumber \\
& \leq & \left( \E \max_{(r-1)D < d \leq rD} G_{d}^{q} H_{d} C^{2(q-1)d} \right)^{2-q} \cdot \left( \sum_{(r-1)D < d \leq rD} \E \frac{G_{d}^{q-1} H_{d}}{C^{2(2-q)d}} \sum_{C^{d-1} \leq |m| \leq C^{d}} G_{d}(m) H_{d}(m) \right)^{q-1} . \nonumber
\end{eqnarray}

We remark that the motivation for the uneven splitting of the Euler products here (moving $G_{d}^{(q-1)^2}$ and $H_{d}^{q-1}$ into the second bracket) is that, as noted above, $\E G_{d}^{q-1} G_{d}(m)$ will behave in approximately the same way as $\E G_{d}^q$, but on the large primes the shifts $im/\mathcal{X}$ provide extra cancellation in $\E H_d H_{d}(m)$. So the best way to split up the $G_d$ terms is ``evenly'', i.e. such that the total exponent of $G_d$ terms in both brackets after H\"{o}lder's inequality remains $q$, whereas for the $H_d$ terms it is better to move a larger piece inside the second bracket (with the sum over $m$) to maximise the cancellation we pick up. As we shall see, the powers of $C$ that we have introduced will serve to balance the final sizes of all the terms.

Using the independence of $f(p)$ for different primes, together with Euler Product Result 1, the sums inside the second bracket are
\begin{eqnarray}
&& \sum_{(r-1)D < d \leq rD} \frac{1}{C^{2(2-q)d}} \sum_{C^{d-1} \leq |m| \leq C^{d}} \E G_{d}^{q-1}G_{d}(m) \E H_{d} H_{d}(m) \nonumber \\
& \ll & \sum_{(r-1)D < d \leq rD} \frac{1}{C^{2(2-q)d}} \sum_{C^{d-1} \leq |m| \leq C^{d}} (1 + \frac{\mathcal{X}}{C^d})^{1+(q-1)^2} (\frac{\mathcal{X}}{1+|m|})^{2(q-1)} (\min\{\mathcal{X},C^d\})^2 (\frac{\min\{\mathcal{X},C^d\}}{1+|m|})^2 \nonumber \\
& \ll & \sum_{(r-1)D < d \leq rD} \frac{C^{O(1)}}{C^{2(2-q)d}} \mathcal{X}^{q^2} C^{(3-q^2)d} = C^{O(1)} \mathcal{X}^{q^2} \sum_{(r-1)D < d \leq rD} C^{-(q-1)^2 d} = C^{O(1)} \mathcal{X}^{q^2} C^{-(q-1)^{2}(r-1)D} .  \nonumber
\end{eqnarray}
When performing this calculation, we noted that the contribution to $\E H_{d} H_{d}(m)$ from primes $p > e^{\mathcal{X}}$ is uniformly bounded (by the first part of Euler Product Result 1), similarly as in our analysis of the case $q \geq 2$. Some of our estimates here were a bit crude, but there seems to be no way to avoid losing some factors $C^{O(1)}$, which further explains why our choice of $C = e^{1/(q-1)}$ is essentially the largest we can make without incurring unacceptable losses.

To bound $\E \max_{(r-1)D < d \leq rD} G_{d}^{q} H_{d} C^{2(q-1)d}$, where we need to handle the maximum in a non-trivial way rather than replacing it by a sum (because this term will not be raised to the small power $q-1$), we will use Probability Result 3. To do this, we first rewrite $\E \max_{(r-1)D < d \leq rD} G_{d}^{q} H_{d} C^{2(q-1)d}$ as
$$ \E |F_{k}(\frac{1}{2} + \sigma)|^2 \max_{(r-1)D < d \leq rD} (G_{d} C^{2d})^{q-1} = \E |F_k(\frac{1}{2}+\sigma)|^{2} \cdot \tilde{\E} \max_{(r-1)D < d \leq rD} \prod_{p \leq e^{\frac{\mathcal{X}}{C^d}}}  \left|1 - \frac{f(p)}{p^{\frac{1}{2}+\sigma}}\right|^{-2(q-1)} C^{2(q-1)d} , $$
where $\tilde{\E}$ is expectation under the ``tilted'' measure defined by $\tilde{\p}(A) = \frac{\E |F_k(1/2+\sigma)|^{2} \textbf{1}_{A}}{\E |F_k(1/2+\sigma)|^{2}}$ for each event $A$ (and $\textbf{1}$ denotes the indicator function). We note, for use in a little while, that if $A$ is an event not involving certain primes then those terms factor out from the expectation and cancel between numerator and denominator in the definition of $\tilde{\p}(A)$. Furthermore, the random variables $f(p)$ are still independent under the measure $\tilde{\p}$, since if $A,B$ are events involving disjoint sets of primes then we can split up the Euler product $|F_k(1/2+\sigma)|^{2}$ into sub-products over the corresponding sets, and then the expectation $\E$ will split up correspondingly.

Now Euler Product Result 1 implies that $\E |F_k(\frac{1}{2}+\sigma)|^{2} \asymp \mathcal{X}$. Furthermore, if we write $L_{d} := \prod_{p \leq e^{\mathcal{X}/C^d}}  \left|1 - \frac{f(p)}{p^{1/2+\sigma}}\right|^{-2(q-1)/1.01}$ (say), and $\lambda_{d} := \tilde{\E} L_d$, then Euler Product Result 1 and the independence of the $f(p)$ imply that
$$ \lambda_{d} = \frac{\E \prod_{p \leq e^{\mathcal{X}/C^d}}  \left|1 - \frac{f(p)}{p^{1/2+\sigma}}\right|^{-2(1+(q-1)/1.01)}}{\E \prod_{p \leq e^{\mathcal{X}/C^d}}  \left|1 - \frac{f(p)}{p^{1/2+\sigma}}\right|^{-2}} \asymp (1 + \frac{\mathcal{X}}{C^d})^{\frac{2(q-1)}{1.01} + \frac{(q-1)^2}{1.01^2}} . $$
We similarly get that $\tilde{\E} L_{d}^{1.01} \asymp (1 + \frac{\mathcal{X}}{C^d})^{2(q-1) + (q-1)^2}$. So we have shown that
\begin{eqnarray}
\E \max_{(r-1)D < d \leq rD} G_{d}^{q} H_{d} C^{2(q-1)d} & \asymp & \mathcal{X} \cdot \tilde{\E} \max_{(r-1)D < d \leq rD} \left(\frac{L_d}{\lambda_{d}}\right)^{1.01} (C^d + \mathcal{X})^{2(q-1)} (1 + \frac{\mathcal{X}}{C^d})^{\frac{(q-1)^2}{1.01}} \nonumber \\
& \ll & \mathcal{X}^{1+2(q-1)} (1 + \frac{\mathcal{X}}{C^{(r-1)D + 1}})^{\frac{(q-1)^2}{1.01}} \tilde{\E} \max_{(r-1)D < d \leq rD} \left(\frac{L_d}{\lambda_{d}}\right)^{1.01} , \nonumber
\end{eqnarray}
where we used the fact that $C^d \leq C \mathcal{X}$ (given our convention that those $d$ for which $C^{d-1} > \mathcal{X}$ are omitted) and $C^{2(q-1)} \ll 1$. Finally, since the $f(p)$ are independent under the measure $\tilde{\p}$ (and so the ``increments'' of different primes in the Euler product are independent), the sequence of random variables $\left(\frac{L_{rD}}{\lambda_{rD}}\right), \left(\frac{L_{rD-1}}{\lambda_{rD-1}}\right), ..., \left(\frac{L_{(r-1)D+1}}{\lambda_{(r-1)D+1}}\right)$ (taken in that order) form a non-negative submartingale relative to $\tilde{\p}$ and to the sigma algebras generated by $(f(p))_{p \leq e^{\mathcal{X}/C^{rD}}}, (f(p))_{p \leq e^{\mathcal{X}/C^{rD-1}}}, ..., (f(p))_{p \leq e^{\mathcal{X}/C^{(r-1)D+1}}}$. Thus Probability Result 3 is applicable, and gives that $\tilde{\E} \max_{(r-1)D < d \leq rD} \left(\frac{L_d}{\lambda_{d}}\right)^{1.01}$ is
$$ \ll \tilde{\E} \left(\frac{L_{(r-1)D + 1}}{\lambda_{(r-1)D + 1}}\right)^{1.01} \asymp \frac{(1 + \frac{\mathcal{X}}{C^{(r-1)D+1}})^{2(q-1) + (q-1)^2}}{\lambda_{(r-1)D + 1}^{1.01}} \asymp (1 + \frac{\mathcal{X}}{C^{(r-1)D+1}})^{(q-1)^2 - \frac{(q-1)^2}{1.01}} . $$

Putting together \eqref{upperq-1}, \eqref{maxdsumm}, and the above calculations, we get that $\vert\vert \int_{-1/2}^{1/2} |F_k(1/2+\sigma + it)|^2 dt \vert\vert_{q}$ is
$$ \ll \frac{1}{\mathcal{X}} \left( \mathcal{X} \; D^{q-1} \sum_{r \leq \frac{(q-1)\log\mathcal{X}}{D} + 1} \left(\frac{\mathcal{X}^{q^2}}{C^{(r-1)D(q-1)^2}} \right)^{2-q} \left(C^{O(1)} \mathcal{X}^{q^2} C^{-(q-1)^{2}(r-1)D} \right)^{q-1} \right)^{1/q} . $$
Recalling that $C = e^{1/(q-1)}$ and collecting terms together, we find this is all
$$ \ll \frac{1}{\mathcal{X}} ( \mathcal{X}^{1+q^2} \; D^{q-1} \sum_{r \leq \frac{(q-1)\log\mathcal{X}}{D} + 1} C^{-(r-1)D(q-1)^2} )^{1/q} = \frac{1}{\mathcal{X}} ( \mathcal{X}^{1+q^2} \; D^{q-1} \sum_{r \leq \frac{(q-1)\log\mathcal{X}}{D} + 1} e^{-(r-1)D(q-1)} )^{1/q} . $$
So if we finally choose $D := \lfloor \frac{1}{q-1} \rfloor$, then both the sum over $r$ and the term $D^{q-1}$ will be $\ll 1$. Recalling that $\mathcal{X} = \min\{\frac{\log x}{e^{k+1}}, \frac{1}{|\sigma|}\}$, we see this bound is as claimed in Key Proposition 1.
\qed

\subsection{Proof of Key Proposition 2, for $1 < q < 2$}
We again begin with the second part of the proposition, where $|N| \geq 1$. In this case we can analyse the terms $\E |F_k(\frac{1}{2}+\sigma + i(\frac{n}{\mathcal{X}} + N + t))|^2 \left( \sum_{|m| \leq \frac{\mathcal{X}}{2} + 1} |F_k(\frac{1}{2}+\sigma + i(\frac{m}{\mathcal{X}} + N + t))|^2 \right)^{q-1}$ in \eqref{radupperq-1} by splitting the sum over $m$ into subsums where $C^{d-1} \leq |m-n| \leq C^d$, and otherwise following the argument from the Steinhaus case. We obtain the same estimates as there, except the error term in Euler Product Result 2 produces an additional factor
$$ \min\{1 + \frac{\mathcal{X}}{C^d}, |N|^{1/100}\}^{|(q-1)^2 - (q-1)| + 4(q-1)} \min\{C^d, \mathcal{X}, 1 + \frac{|N|^{1/100}}{1 + \mathcal{X}/C^d} \}^{4} \ll \min\{\mathcal{X}, |N|^{1/100}\}^4 $$
when estimating (the analogue of) the terms $\E G_{d}^{q-1}G_{d}(m) \E H_{d} H_{d}(m)$, and an additional factor $\min\{1 + \frac{\mathcal{X}}{C^{(r-1)D + 1}}, |N|^{1/100}\}^{q(q-1)} \ll \min\{\mathcal{X}, |N|^{1/100}\}^{q(q-1)}$ when estimating (the analogue of) the term $\E \max_{(r-1)D < d \leq rD} G_{d}^{q} H_{d} C^{2(q-1)d}$. So overall we get the same bound as in the Steinhaus case, apart from a factor
$$ (\min\{\mathcal{X}, |N|^{1/100}\}^{q(q-1)})^{2-q} (\min\{\mathcal{X}, |N|^{1/100}\}^{4})^{(q-1)} \ll \min\{\mathcal{X}, |N|^{1/100}\}^{(q-1)(4 + 2q - q^2)} . $$
A small calculation shows that for $1 \leq q \leq 2$, we have $(q-1)(4 + 2q - q^2) \leq 5(q-1) \leq q(q+1)$, giving the factor $\min\{\mathcal{X}, |N|^{1/100}\}^{q(q+1)} = \min\{\frac{\log x}{e^{k+1}}, \frac{1}{|\sigma|}, |N|^{1/100}\}^{q(q+1)}$ claimed in the second part of Key Proposition 2.

When $N=0$, to prove Key Proposition 2 we need to bound
$$ \sum_{|n| \leq \frac{\mathcal{X}}{2} + 1} \E|F_k(\frac{1}{2}+\sigma + i(\frac{n}{\mathcal{X}} + t))|^2 \left( \sum_{|m| \leq \frac{\mathcal{X}}{2} + 1} |F_k(1/2+\sigma + i(\frac{m}{\mathcal{X}} + t))|^2 \right)^{q-1} $$
in \eqref{radupperq-1}. Following the same argument that led to the bound \eqref{maxdsumm} in the Steinhaus case, but now splitting the sum over $m$ according to the size of $|m|-|n|$ rather than the size of $|m|$, one obtains that
\begin{eqnarray}\label{radmaxdsumm}
&& \sum_{|n| \leq \frac{\mathcal{X}}{2} + 1} \E|F_k(\frac{1}{2}+\sigma + i(\frac{n}{\mathcal{X}} + t))|^2 \left( \sum_{|m| \leq \frac{\mathcal{X}}{2} + 1} |F_k(1/2+\sigma + i(\frac{m}{\mathcal{X}} + t))|^2 \right)^{q-1} \nonumber \\
& \ll & D^{q-1} \sum_{|n| \leq \frac{\mathcal{X}}{2} + 1} \sum_{r \leq \frac{(q-1)\log\mathcal{X}}{D} + 1} \E \max_{(r-1)D < d \leq rD} |F_k(\frac{1}{2}+\sigma + i(\frac{n}{\mathcal{X}} + t))|^2 \cdot \nonumber \\
&& \cdot \Biggl( \sum_{\substack{ |m| \leq \frac{\mathcal{X}}{2} + 1 , \\ C^{d-1} \leq ||m| - |n|| \leq C^{d}}} |F_k(\frac{1}{2}+\sigma + i(\frac{m}{\mathcal{X}} + t))|^2 \Biggr)^{q-1} .
\end{eqnarray}
Here we again have $C = e^{1/(q-1)}$, and $D = D(q) \in \N$ is a parameter, and we adopt our usual conventions (analogously to the Steinhaus case) about including the $|m|=|n|$ term when $d=1$ and omitting overly large terms from all sums.

Now for each $d \geq 1$ and $|m| \leq \frac{\mathcal{X}}{2} + 1$, and treating $t$ as fixed, we shall set 
$$ G_{d}(m) := \prod_{p \leq e^{\mathcal{X}/C^d}}  \left|1 + \frac{f(p)}{p^{1/2+\sigma+i(m/\mathcal{X} + t)}}\right|^{2} , \;\;\;\;\; \text{and} \;\;\;\;\; H_{d}(m) := \frac{|F_k(\frac{1}{2}+\sigma + i(\frac{m}{\mathcal{X}} + t))|^2}{G_{d}(m)} . $$
This is the same notation that we used in the Steinhaus case, but with the Euler products now replaced by their Rademacher versions (supported on squarefree numbers only). Splitting the expectation and applying H\"{o}lder's inequality as in the Steinhaus case, it follows that \eqref{radmaxdsumm} is
\begin{eqnarray}
& \ll & D^{q-1} \sum_{|n| \leq \frac{\mathcal{X}}{2} + 1} \sum_{r \leq \frac{(q-1)\log\mathcal{X}}{D} + 1} \left( \E \max_{(r-1)D < d \leq rD} G_{d}(n)^{q} H_{d}(n) C^{2(q-1)d} \right)^{2-q} \cdot \nonumber \\
&& \cdot \Biggl( \sum_{(r-1)D < d \leq rD} \E \frac{G_{d}(n)^{q-1} H_{d}(n)}{C^{2(2-q)d}} \sum_{\substack{ |m| \leq \frac{\mathcal{X}}{2} + 1 , \\ C^{d-1} \leq ||m| - |n|| \leq C^{d}}} G_{d}(m) H_{d}(m) \Biggr)^{q-1} . \nonumber
\end{eqnarray}

Continuing to follow the argument from the Steinhaus case, but using Euler Product Result 2 in place of Euler Product Result 1, we can bound these terms further. Proceeding to do this, and noting that one of the terms $|m-n|, |m+n|$ that arise in Euler Product Result 2 will always equal $||m|-|n||$ and the other will equal $|m|+|n| \geq |n|$, we find the sums in the second bracket are
\begin{eqnarray}
&& \sum_{(r-1)D < d \leq rD} \frac{1}{C^{2(2-q)d}} \sum_{\substack{ |m| \leq \frac{\mathcal{X}}{2} + 1 , \\ C^{d-1} \leq ||m| - |n|| \leq C^{d}}} \E G_{d}(n)^{q-1}G_{d}(m) \E H_{d}(n) H_{d}(m) \nonumber \\
& \ll & \sum_{(r-1)D < d \leq rD} \frac{1}{C^{2(2-q)d}} \sum_{C^{d-1} \leq ||m| - |n|| \leq C^{d}} \min\{\frac{\mathcal{X}}{C^d}, \frac{\mathcal{X}}{|n|}\}^{(q-1)^2 - (q-1)} (1 + \frac{\mathcal{X}}{C^d})^{1+(q-1)^2} \cdot \nonumber \\
&& \cdot \left( \frac{\mathcal{X}}{1+||m|-|n||} \min\{\frac{\mathcal{X}}{C^d}, \frac{\mathcal{X}}{|n|}\} \right)^{2(q-1)} (\min\{\mathcal{X},C^d\})^2 \left( \frac{\min\{\mathcal{X},C^d\}}{1+||m|-|n||} (1 + \frac{\min\{\mathcal{X},C^d\}}{1+|m|+|n|}) \right)^2 . \nonumber
\end{eqnarray}
Collecting terms together, and then upper bounding $\min\{\frac{\mathcal{X}}{C^d}, \frac{\mathcal{X}}{|n|}\}^{q(q-1)}$ by $(\frac{\mathcal{X}}{1+|n|})^{q(q-1)}$ and upper bounding $\min\{\mathcal{X},C^d\}$ everywhere else by $C^d$, the above is
$$ \ll (\frac{\mathcal{X}}{1+|n|})^{q(q-1)} \sum_{(r-1)D < d \leq rD} \frac{C^{O(1)}}{C^{2(2-q)d}} \mathcal{X}^{q^2} C^{(3-q^2)d} = (\frac{\mathcal{X}}{1+|n|})^{q(q-1)} C^{O(1)} \mathcal{X}^{q^2} C^{-(q-1)^{2}(r-1)D} . $$
We can also adapt the Steinhaus argument to bound $\E \max_{(r-1)D < d \leq rD} G_{d}(n)^{q} H_{d}(n) C^{2(q-1)d}$. In this case we again have $\E |F_k(\frac{1}{2}+\sigma + i(\frac{n}{\mathcal{X}} + t))|^2 \asymp \mathcal{X}$, and we may define the ``tilted'' measure $\tilde{\p}$ and set $L_{d} := \prod_{p \leq e^{\mathcal{X}/C^d}}  \left|1 + \frac{f(p)}{p^{1/2+\sigma+ i(\frac{n}{\mathcal{X}} + t)}}\right|^{2(q-1)/1.01}$ analogously to the Steinhaus case. Then Euler Product Result 2 implies that $\tilde{\E} L_d \asymp \min\{\frac{\mathcal{X}}{C^d}, \frac{\mathcal{X}}{|n|}\}^{(1 + \frac{q-1}{1.01}) \frac{q-1}{1.01}} (1 + \frac{\mathcal{X}}{C^d})^{\frac{2(q-1)}{1.01} + \frac{(q-1)^2}{1.01^2}}$ and $\tilde{\E} L_d^{1.01} \asymp \min\{\frac{\mathcal{X}}{C^d}, \frac{\mathcal{X}}{|n|}\}^{q(q-1)} (1 + \frac{\mathcal{X}}{C^d})^{2(q-1) + (q-1)^2}$. So the same submartingale argument as in the Steinhaus case shows that
\begin{eqnarray}
\E \max_{(r-1)D < d \leq rD} G_{d}(n)^{q} H_{d}(n) C^{2(q-1)d} & \ll & \mathcal{X}^{1+2(q-1)} (1 + \frac{\mathcal{X}}{C^{(r-1)D + 1}})^{(q-1)^2} \min\{\frac{\mathcal{X}}{C^{(r-1)D + 1}}, \frac{\mathcal{X}}{|n|}\}^{q(q-1)} \nonumber \\
& \ll & \frac{\mathcal{X}^{q^2}}{C^{(r-1)D(q-1)^2}} (\frac{\mathcal{X}}{1+|n|})^{q(q-1)} . \nonumber
\end{eqnarray}

Putting everything together, recalling that $C = e^{1/(q-1)}$ and choosing $D := \lfloor \frac{1}{q-1} \rfloor$, we deduce that \eqref{radmaxdsumm} is
$$ \ll D^{q-1} \sum_{|n| \leq \frac{\mathcal{X}}{2} + 1} \sum_{r \leq \frac{(q-1)\log\mathcal{X}}{D} + 1} \frac{\mathcal{X}^{q^2}}{C^{(r-1)D(q-1)^2}} (\frac{\mathcal{X}}{1+|n|})^{q(q-1)} \ll \mathcal{X}^{q^2} \sum_{|n| \leq \frac{\mathcal{X}}{2} + 1} (\frac{\mathcal{X}}{1+|n|})^{q(q-1)} . $$
Since $q_0 = (1 + \sqrt{5})/2 \approx 1.618$ satisfies $q_0 (q_0 - 1) = 1$, and we have $1 < q < 2$, the sum over $n$ here is $\ll \mathcal{X}^{\max\{1, q(q-1)\}} \min\{\log\mathcal{X}, \frac{1}{|q-q_0|}\}$. Substituting into \eqref{radupperq-1}, and recalling that $\mathcal{X} := \min\{\frac{\log x}{e^{k+1}}, \frac{1}{|\sigma|}\}$, this gives the first ($N=0$) bound claimed in Key Proposition 2 when $1 < q < 2$.
\qed

\section{Proofs of the lower bounds in Theorems 1 and 2}\label{secmainlower}
Recall that $F(s)$ denotes the Euler product of $f(n)$ over $x$-smooth numbers.

\subsection{The lower bound in the Steinhaus case}\label{rmfiisteinhauslowermain}
To prove the lower bound part of Theorem 1, in view of Proposition \ref{propstlower} our main work will be to prove a suitable lower bound for $\vert\vert \int_{-1/2}^{1/2} |F(1/2 + \frac{4Vq}{\log x} + it)|^2 dt \vert\vert_{q}$, where $V$ is a large fixed constant. (We also need an upper bound for the subtracted quantity $\vert\vert \int_{-1/2}^{1/2} |F(1/2 + \frac{2Vq}{\log x} + it)|^2 dt \vert\vert_{q}$, but this will follow directly from Key Proposition 1.)

To obtain our lower bound, we note first that for any $q \geq 1$ we have
\begin{eqnarray}
\left( \int_{- 1/2}^{1/2} |F(1/2+ \frac{4Vq}{\log x} + it)|^2 dt \right)^q & \geq & \left( \sum_{|k| \leq \frac{\log x - 1}{2}} \int_{\frac{k-1/2}{\log x}}^{\frac{k+1/2}{\log x}} |F(1/2+ \frac{4Vq}{\log x} + it)|^2 dt \right)^q \nonumber \\
& \geq & \sum_{|k| \leq \frac{\log x - 1}{2}} \left( \int_{\frac{k-1/2}{\log x}}^{\frac{k+1/2}{\log x}} |F(1/2+ \frac{4Vq}{\log x} + it)|^2 dt \right)^q . \nonumber
\end{eqnarray}
This step could be wasteful if many of the pieces $\int_{\frac{k-1/2}{\log x}}^{\frac{k+1/2}{\log x}} |F(1/2+ \frac{4Vq}{\log x} + it)|^2 dt$ made substantial contributions to the full integral. But for large $q$ we expect instead that the dominant contribution should come from just a few large (and therefore rare) contributions, so we will not lose too much. In the Steinhaus case, since the distribution of $\int_{\frac{k-1/2}{\log x}}^{\frac{k+1/2}{\log x}} |F(1/2+ \frac{4Vq}{\log x} + it)|^2 dt$ is independent of $k$ we get the simpler lower bound
$$ \E\left( \int_{- 1/2}^{1/2} |F(1/2+ \frac{4Vq}{\log x} + it)|^2 dt \right)^q \gg \log x \cdot \E\left( \int_{-\frac{1}{2\log x}}^{\frac{1}{2\log x}} |F(1/2+ \frac{4Vq}{\log x} + it)|^2 dt \right)^q . $$

Now we want to remove the remaining short integral over $t$, which is a technical obstacle to connecting up the expectation with the random product $F$. Heuristically, since the Euler product shouldn't vary much on intervals of length $1/\log x$ we should simply obtain something like $\left( \frac{1}{\log x} |F(1/2 + \frac{4Vq}{\log x})|^2 \right)^q$ in the bracket. It turns out that a neat way to handle this issue is using Jensen's inequality (applied to the normalised integral $\int_{-\frac{1}{2\log x}}^{\frac{1}{2\log x}} \log x \; dt$), which implies that $\E \left( \int_{- \frac{1}{2\log x}}^{\frac{1}{2\log x}} |F(1/2+\frac{4Vq}{\log x} + it)|^2 dt \right)^q$ is
\begin{eqnarray}
& = & \frac{1}{\log^{q}x} \E \left( \int_{- \frac{1}{2\log x}}^{\frac{1}{2\log x}} \log x \cdot e^{2\log|F(1/2+4Vq/\log x + it)|} dt \right)^q \nonumber \\
& \geq & \frac{1}{\log^{q}x} \E \left( \exp\{\int_{- \frac{1}{2\log x}}^{\frac{1}{2\log x}} \log x \cdot 2\log|F(1/2+\frac{4Vq}{\log x} + it)| dt\} \right)^q \nonumber \\
& = & \frac{1}{\log^{q}x} \E \exp\{2q \int_{- \frac{1}{2\log x}}^{\frac{1}{2\log x}} \log x \cdot \log|F(1/2+\frac{4Vq}{\log x} + it)| dt \} . \nonumber
\end{eqnarray}
Here the exponential, inside the expectation, may be rewritten as
\begin{eqnarray}
&& \prod_{p \leq x} \exp\Biggl\{-2q \int_{- \frac{1}{2\log x}}^{\frac{1}{2\log x}} \log x \cdot \Re\log\Biggl(1 - \frac{f(p)}{p^{1/2+4Vq/\log x + it}}\Biggr) dt \Biggr\} \nonumber \\
& = & \prod_{p \leq x} \exp\Biggl\{2q \log x \Re\Biggl( \frac{f(p)}{p^{1/2+\frac{4Vq}{\log x}}} \int_{- \frac{1}{2\log x}}^{\frac{1}{2\log x}}e^{-it\log p} dt + \frac{f(p)^2}{2 p^{1+\frac{8Vq}{\log x}}} \int_{- \frac{1}{2\log x}}^{\frac{1}{2\log x}}e^{-2it\log p} dt \Biggr) + O(\frac{q}{p^{3/2}}) \Biggr\} \nonumber \\
& = & e^{O(q)} \prod_{p \leq x} \exp\Biggl\{2q \Re \Biggl( \frac{f(p)}{p^{1/2+4Vq/\log x}} \log x \int_{- \frac{1}{2\log x}}^{\frac{1}{2\log x}}e^{-it\log p} dt + \frac{f(p)^2}{2 p^{1+8Vq/\log x}} \Biggr)\Biggr\} . \nonumber
\end{eqnarray}
The first equality here uses the Taylor expansion of the logarithm, and the second equality uses the estimate $\int_{- \frac{1}{2\log x}}^{\frac{1}{2\log x}}e^{-2it\log p} dt = \int_{- \frac{1}{2\log x}}^{\frac{1}{2\log x}} (1 + O(|t|\log p)) dt = 1/\log x + O((\log p)/\log^{2}x)$ (and also the fact that $\sum_{p \leq x} \frac{\log p}{p^{1+8Vq/\log x}} \ll \log x$). Putting things together, using the independence of $f(p)$ for different primes $p$ to move the expectation inside the product, we have shown that our original object $\E\left( \int_{- 1/2}^{1/2} |F(1/2+\frac{4Vq}{\log x} + it)|^2 dt \right)^q$ is
\begin{equation}\label{stlowerstep1}
\geq \frac{e^{O(q)}}{\log^{q-1}x} \prod_{p \leq x} \E \exp\Biggl\{2q \Re \Biggl( \frac{f(p)}{p^{1/2+4Vq/\log x}} \log x \int_{- \frac{1}{2\log x}}^{\frac{1}{2\log x}}e^{-it\log p} dt + \frac{f(p)^2}{2 p^{1+8Vq/\log x}} \Biggr)\Biggr\} .
\end{equation}

It will be convenient to note some simple bounds for the quantity inside the exponential, which we will use shortly. Firstly, this quantity is always trivially bounded by $O(q/\sqrt{p})$. Secondly, using our previous calculation that $\int_{- \frac{1}{2\log x}}^{\frac{1}{2\log x}}e^{-it\log p} dt = 1/\log x + O((\log p)/\log^{2}x)$, we can obtain that
$$ 2q \Re \Biggl( \frac{f(p)}{p^{1/2+\frac{4Vq}{\log x}}} \log x \int_{- \frac{1}{2\log x}}^{\frac{1}{2\log x}}e^{-it\log p} dt + \frac{f(p)^2}{2 p^{1+\frac{8Vq}{\log x}}} \Biggr) = \frac{2q \Re f(p)}{p^{1/2+\frac{4Vq}{\log x}}} + O\Biggl(\frac{q\log p}{p^{1/2+\frac{4Vq}{\log x}} \log x} + \frac{q}{p} \Biggr) . $$

To conclude, we note that certainly when $100q^2 \leq p \leq x$ we have, in view of the Taylor expansion of the exponential and the simple bounds noted above and the fact that $\E (\Re f(p))^2 = 1/2$, that
\begin{eqnarray}
&& \E \exp\Biggl\{2q \Re \Biggl( \frac{f(p)}{p^{1/2+4Vq/\log x}} \log x \int_{- \frac{1}{2\log x}}^{\frac{1}{2\log x}}e^{-it\log p} dt + \frac{f(p)^2}{2 p^{1+8Vq/\log x}} \Biggr)\Biggr\} \nonumber \\
& = & \E \left( 1 + 2q \Re \Biggl( \frac{f(p)}{p^{1/2+\frac{4Vq}{\log x}}} \log x \int_{- \frac{1}{2\log x}}^{\frac{1}{2\log x}}e^{-it\log p} dt + \frac{f(p)^2}{2 p^{1+\frac{8Vq}{\log x}}} \Biggr) + \frac{2q^2 (\Re f(p))^2}{p^{1 + \frac{8Vq}{\log x}}} + \right. \nonumber \\
&& \left. + O(\frac{q^2 \log p}{p^{1 + \frac{8Vq}{\log x}} \log x} + \frac{q^3}{p^{3/2}}) \right) \nonumber \\
& = & 1 + \frac{q^2}{p^{1+8Vq/\log x}} + O\Biggl(\frac{q^2 \log p}{p^{1 + 8Vq/\log x} \log x} + \frac{q^3}{p^{3/2}}\Biggr) . \nonumber
\end{eqnarray}
When $p < 100q^2$, we shall instead use the trivial bound $\exp\left\{ O(\frac{q}{\sqrt{p}}) \right\}$. Inserting these into \eqref{stlowerstep1}, we get an overall lower bound
\begin{eqnarray}
&& \E\left( \int_{- 1/2}^{1/2} |F(1/2+\frac{4Vq}{\log x} + it)|^2 dt \right)^q \nonumber \\
& \geq & \frac{e^{O(q)}}{\log^{q-1}x} \prod_{p < 100q^2} \exp\left\{ O(\frac{q}{\sqrt{p}}) \right\} \prod_{100q^2 \leq p \leq x} \exp\left\{\frac{q^2}{p^{1+8Vq/\log x}} + O(\frac{q^2 \log p}{p^{1 + 8Vq/\log x} \log x} + \frac{q^3}{p^{3/2}}) \right\} \nonumber \\
& = & \frac{e^{O(q^{2}/\log(2q))}}{\log^{q-1}x} \prod_{100q^2 \leq p \leq x} \exp\left\{\frac{q^2}{p^{1+8Vq/\log x}} \right\} = \frac{e^{O(q^{2})}}{\log^{q-1}x} \left(\frac{\log x}{Vq \log(2q)}\right)^{q^2} . \nonumber
\end{eqnarray}
Here the final equality follows because, similarly as in the calculations in section \ref{seceasycases}, we have $\prod_{100q^2 \leq p \leq x} \exp\left\{\frac{q^2}{p^{1+8Vq/\log x}} \right\} = e^{O(q^2)} \zeta(1 + \frac{8Vq}{\log x})^{q^2} e^{-\sum_{p < 100q^2} \frac{q^2}{p^{1+8Vq/\log x}}}$. Then $\zeta(1 + \frac{8Vq}{\log x})^{q^2} = e^{O(q^2)} (\frac{\log x}{Vq})^{q^2}$, and $e^{-\sum_{p < 100q^2} \frac{q^2}{p^{1+8Vq/\log x}}} = \frac{e^{O(q^2)}}{\log^{q^2}(2q)}$ by Mertens' estimates for sums over primes.

Inserting this into Proposition \ref{propstlower} we find that $\vert\vert \sum_{n \leq x} f(n) \vert\vert_{2q}$ is
$$ \gg \sqrt{\frac{x}{\log x}} \Biggl(\frac{e^{O(q)}}{\log^{(q-1)/2q}x} \left(\frac{\log x}{Vq \log(2q)}\right)^{q/2}  - \frac{C}{e^{Vq/2}} \vert\vert \int_{-1/2}^{1/2} |F(1/2 + \frac{2Vq}{\log x} + it)|^2 dt \vert\vert_{q}^{1/2} \Biggr) . $$
And using Key Proposition 1 with $k=-1$ to control the subtracted term, provided that $\frac{2Vq}{\log x} \leq \frac{1}{100\log(2q)}$ we can lower bound everything by
$$ \sqrt{\frac{x}{\log x}} \Biggl(\frac{e^{O(q)}}{\log^{(q-1)/2q}x} \left(\frac{\log x}{Vq \log(2q)}\right)^{q/2}  - \frac{C e^{O(q)}}{e^{Vq/2}} \frac{(Vq)^{(q-1)/2q}}{\log^{(q-1)/2q}x} \left(\frac{\log x}{Vq \log(2q)}\right)^{q/2} \Biggr) . $$
If we set $V$ to be a sufficiently large fixed constant, the subtracted term will be negligible compared with the first term, and our Theorem 1 lower bound will be proved. It only remains to note that the condition $\frac{2Vq}{\log x} \leq \frac{1}{100\log(2q)}$ is then satisfied provided $q \leq \frac{c\log x}{\log\log x}$, for a sufficiently small fixed constant $c > 0$. 
\qed

\subsection{The lower bound in the Rademacher case}
To prove the lower bound part of Theorem 2, we shall invoke Proposition \ref{propradlower} and adapt the argument from the previous subsection to lower bound $\vert\vert \int_{-1/2}^{1/2} |F(1/2 + \frac{4Vq}{\log x} + it)|^2 dt \vert\vert_{q}$, where $V$ is a large fixed constant and now $F(s)$ denotes the Rademacher random Euler product.

Indeed, exactly the same argument as in section \ref{rmfiisteinhauslowermain} gives, for any $q \geq 1$, that
$$ \E \left( \int_{- 1/2}^{1/2} |F(1/2+ \frac{4Vq}{\log x} + it)|^2 dt \right)^q \geq \sum_{|k| \leq \frac{\log x - 1}{2}} \E \left( \int_{\frac{k-1/2}{\log x}}^{\frac{k+1/2}{\log x}} |F(1/2+ \frac{4Vq}{\log x} + it)|^2 dt \right)^q . $$
Using Jensen's inequality, also as in section \ref{rmfiisteinhauslowermain}, shows this is all
$$ \geq \frac{1}{\log^{q}x} \sum_{|k| \leq \frac{\log x - 1}{2}} \E \exp\{2q \int_{\frac{k - 1/2}{\log x}}^{\frac{k + 1/2}{\log x}} \log x \cdot \log|F(1/2+\frac{4Vq}{\log x} + it)| dt \} . $$
And recalling that in the Rademacher case we have $F(s) = \prod_{p \leq x} (1 + \frac{f(p)}{p^s})$, with $f(p) \in \{\pm 1\}$, we find this is all
\begin{eqnarray}
& = & \frac{1}{\log^{q}x} \sum_{|k| \leq \frac{\log x - 1}{2}} \prod_{p \leq x} \E \exp\Biggl\{2q \int_{\frac{k - 1/2}{\log x}}^{\frac{k + 1/2}{\log x}} \log x \cdot \Re \log\Biggl( 1 + \frac{f(p)}{p^{1/2+4Vq/\log x + it}} \Biggr) dt \Biggr\} \nonumber \\
& = & \frac{1}{\log^{q}x} \sum_{|k| \leq \frac{\log x - 1}{2}} \prod_{p \leq x} \E \exp\Biggl\{2q \int_{\frac{k - 1/2}{\log x}}^{\frac{k + 1/2}{\log x}} \log x \cdot \Biggl( \frac{f(p) \cos(t\log p)}{p^{1/2+4Vq/\log x}} - \frac{\cos(2t\log p)}{2p^{1+8Vq/\log x}} + O(\frac{1}{p^{3/2}}) \Biggr) dt \Biggr\} . \nonumber
\end{eqnarray}
At this stage we cannot efficiently remove the integral of $\cos(t\log p)$ in the first term, but for the second term we can write $\cos(2t\log p) = \cos(\frac{2k\log p}{\log x}) + O(\frac{\log p}{\log x})$. The total contribution from these ``big Oh'' terms for all $p$, as well as from the $O(1/p^{3/2})$ term, is a multiplicative factor $e^{O(q)}$. So we obtain that $\E\left( \int_{- 1/2}^{1/2} |F(1/2+\frac{4Vq}{\log x} + it)|^2 dt \right)^q$ is
\begin{equation}\label{radlowerstep1}
\geq \frac{e^{O(q)}}{\log^{q}x} \sum_{|k| \leq \frac{\log x - 1}{2}} \prod_{p \leq x} \E \exp\Biggl\{2q \Biggl( \frac{f(p)}{p^{1/2+4Vq/\log x}} \log x \int_{\frac{k - 1/2}{\log x}}^{\frac{k + 1/2}{\log x}} \cos(t\log p) dt - \frac{\cos(\frac{2k\log p}{\log x})}{2 p^{1+8Vq/\log x}} \Biggr)\Biggr\} ,
\end{equation}
which is the Rademacher analogue of \eqref{stlowerstep1} from the Steinhaus case.

Next, when $100q^2 \leq p \leq x$ we have, in view of the Taylor expansion of the exponential (and the fact that $f(p)^2 \equiv 1$), that
\begin{eqnarray}
&& \E \exp\Biggl\{2q \Biggl( \frac{f(p)}{p^{1/2+4Vq/\log x}} \log x \int_{\frac{k - 1/2}{\log x}}^{\frac{k + 1/2}{\log x}} \cos(t\log p) dt - \frac{\cos(\frac{2k\log p}{\log x})}{2 p^{1+8Vq/\log x}} \Biggr)\Biggr\} \nonumber \\
& = & \E \left( 1 + 2q \Biggl( \frac{f(p)}{p^{1/2+\frac{4Vq}{\log x}}} \log x \int_{\frac{k - 1/2}{\log x}}^{\frac{k + 1/2}{\log x}} \cos(t\log p) dt - \frac{\cos(\frac{2k\log p}{\log x})}{2 p^{1+\frac{8Vq}{\log x}}} \Biggr) + \frac{2q^2 \cos^2(\frac{k\log p}{\log x})}{p^{1 + \frac{8Vq}{\log x}}} + \right. \nonumber \\
&& \left. + O(\frac{q^2 \log p}{p^{1 + \frac{8Vq}{\log x}} \log x} + \frac{q^3}{p^{3/2}}) \right) . \nonumber
\end{eqnarray}
Using the cosine identity $\cos^2(\frac{k\log p}{\log x}) = (1/2)(1 + \cos(\frac{2k\log p}{\log x}))$, and the fact that $\E f(p) = 0$, we find the above is
$$ = 1 + \frac{q^2 + (q^2 - q)\cos(\frac{2k\log p}{\log x})}{p^{1+8Vq/\log x}} + O\Biggl(\frac{q^2 \log p}{p^{1 + 8Vq/\log x} \log x} + \frac{q^3}{p^{3/2}}\Biggr) . $$
When $p < 100q^2$, we shall instead use the trivial bound $\exp\left\{ O(\frac{q}{\sqrt{p}}) \right\}$. Inserting these into \eqref{radlowerstep1}, we get
$$ \E\left( \int_{- 1/2}^{1/2} |F(1/2+\frac{4Vq}{\log x} + it)|^2 dt \right)^q \geq \frac{e^{O(\frac{q^{2}}{\log(2q)})}}{\log^{q}x} \sum_{|k| \leq \frac{\log x - 1}{2}} \prod_{100q^2 \leq p \leq x} \exp\{ \frac{q^2 + (q^2 - q)\cos(\frac{2k\log p}{\log x})}{p^{1+\frac{8Vq}{\log x}}} \} . $$

Now when $2 \leq q \leq \frac{c\log x}{\log\log x}$, say, we can afford to discard all the terms in this lower bound except the $k=0$ term, which gives us that $ \E\left( \int_{- 1/2}^{1/2} |F(1/2+\frac{4Vq}{\log x} + it)|^2 dt \right)^q$ is $\geq \frac{e^{O(\frac{q^{2}}{\log(2q)})}}{\log^{q}x} \prod_{100q^2 \leq p \leq x} \exp\{ \frac{2q^2 - q}{p^{1+\frac{8Vq}{\log x}}} \} = \frac{e^{O(q^{2})}}{\log^{q}x} \left(\frac{\log x}{Vq \log(2q)}\right)^{2q^2 - q}$. Inserting this into Proposition \ref{propradlower}, and applying Key Proposition 2 with $k=-1$ and $\sigma = 2Vq/\log x$ to control the subtracted term there, we find that $\vert\vert \sum_{n \leq x} f(n) \vert\vert_{2q}$ is
$$ \gg \sqrt{\frac{x}{\log x}} \Biggl(\frac{e^{O(q)}}{\log^{1/2}x} \left(\frac{\log x}{Vq \log(2q)}\right)^{q - 1/2}  - \frac{C e^{O(q)}}{e^{Vq/2}} \frac{(Vq)^{1/2}}{\log^{1/2}x} \left(\frac{\log x}{Vq \log(2q)}\right)^{q - 1/2} \Biggr) . $$
If $V$ is a sufficiently large constant, the subtracted term is negligible compared with the first term and we obtain the lower bound claimed in Theorem 2.

When $1 \leq q < 2$ (or really when $1 \leq q \leq q_0 = (1+\sqrt{5})/2$), we cannot afford to take quite such a crude approach. Using Chebychev's estimates and the Prime Number Theorem as in the proof of Euler Product Result 1, we have
\begin{eqnarray}
\sum_{100q^2 \leq p \leq x} \frac{\cos(\frac{2k\log p}{\log x})}{p^{1+\frac{8Vq}{\log x}}} & = & \sum_{2 \leq p \leq x^{1/V}} \frac{\cos(\frac{2k\log p}{\log x})}{p} + O(1) = \int_{\log 2}^{\log(x^{1/V})} \frac{\cos(\frac{2k}{\log x} u)}{u} du + O(1) \nonumber \\
& = & \log\min\{ \log(x^{1/V}) , \log(x^{1/(1+|k|)}) \} + O(1) . \nonumber
\end{eqnarray}
Using this estimate, we get a lower bound for $\E\left( \int_{- 1/2}^{1/2} |F(1/2+\frac{4Vq}{\log x} + it)|^2 dt \right)^q$ that is $\gg \frac{1}{\log^{q}x} \left(\frac{\log x}{V}\right)^{q^2} \sum_{|k| \leq \frac{\log x - 1}{2}} \min\{ \frac{\log x}{V}, \frac{\log x}{1+|k|} \}^{q^2 - q}$, and (remembering that $q_0$ satisfies $q_{0}^2 - q_0 = 1$) this is $\gg \frac{1}{\log^{q}x} \left(\frac{\log x}{V}\right)^{q^2 + \max\{1,q^{2}-q\}} \min\{\log\log x, \frac{1}{|q-q_0|}\}$. Again, inserting this in Key Proposition 2 produces the lower bound claimed in Theorem 2.
\qed


\end{document}